\DocumentMetadata{testphase=new-or-1}
\documentclass[11pt,a4paper]{article}
\usepackage{lmodern}
\usepackage[margin=1in]{geometry}
\usepackage[authoryear,round,semicolon]{natbib}
\usepackage[utf8]{inputenc}
\usepackage[T1]{fontenc}
\usepackage{microtype}
\usepackage{needspace}
\usepackage{amsfonts}
\usepackage{amsmath, amssymb, amsthm}
\usepackage{mathtools}
\usepackage{array, booktabs}
\usepackage[table]{xcolor}
\usepackage{float}
\usepackage{placeins}
\usepackage{graphicx}
\usepackage{grffile}
\usepackage{hyperref}
\usepackage[nameinlink,capitalise]{cleveref}
\usepackage{listings}
\usepackage{longtable}
\usepackage{makecell}
\usepackage{multirow}
\usepackage{multicol}
\usepackage{pifont}
\usepackage{siunitx}
\usepackage{subcaption}
\usepackage{tabularx}
\usepackage{tikz}

\usetikzlibrary{positioning, arrows.meta, shapes.geometric, fit, calc}
\usepackage{pgfplots}
\usepgfplotslibrary{groupplots}
\pgfplotsset{
  compat=1.16,
  set layers=standard,
  layers/standard/.append style={
    axis line style={/pgfplots/on layer=axis foreground},
    tick style={/pgfplots/on layer=axis foreground},
  },
  legend image post style={mark size=2pt,scale=1.3,opacity=1},
}
\usepackage{standalone}
\usepackage[disable]{todonotes}
\usepackage{algorithm}
\usepackage{algpseudocode}
\usepackage{enumitem}
\setlist{itemsep=2pt, topsep=3pt, parsep=0pt}

\let\cite\citep
\crefname{observation}{Observation}{Observations}
\Crefname{observation}{Observation}{Observations}
\crefname{corollary}{Corollary}{Corollaries}
\Crefname{corollary}{Corollary}{Corollaries}

\usetikzlibrary{shapes.geometric, positioning, calc, fit, patterns}

\definecolor{ibmBlue}{HTML}{648FFF}
\definecolor{ibmPurple}{HTML}{785EF0}
\definecolor{ibmMagenta}{HTML}{DC267F}
\definecolor{ibmOrange}{HTML}{FE6100}
\definecolor{ibmGold}{HTML}{FFB000}

\colorlet{pigMethod0}{ibmBlue}
\colorlet{pigMethod1}{ibmPurple}
\colorlet{pigMethod2}{ibmMagenta}
\colorlet{pigMethod3}{ibmOrange}
\colorlet{pigMethod4}{ibmGold}
\colorlet{pigMethod5}{black}

\colorlet{pigSrc1}{ibmBlue}
\colorlet{pigSrc2}{ibmOrange}
\colorlet{pigSrc3}{ibmGold}
\colorlet{pigSrc4}{ibmPurple}
\colorlet{pigSrc5}{ibmMagenta}
\colorlet{pigDefault}{pigSrc1}
\definecolor{pigNeutral}{HTML}{666666}
\definecolor{pigMerge}{HTML}{7F7F7F}
\colorlet{pigNew}{black}

\hypersetup{linkbordercolor=ibmPurple,citebordercolor=ibmBlue,urlbordercolor=ibmBlue,filebordercolor=ibmBlue}

\newcommand{\pigSetRowSpacing}{%
  \def\pigBarWidth{9pt}%
  \def\pigBarGap{2.5pt}%
  \pgfmathsetlengthmacro{\pigRowPitch}{\pigBarWidth+\pigBarGap}%
  \pgfmathsetlengthmacro{\pigMethodAxisHeight}{5.28*\pigRowPitch}%
  \pgfmathsetlengthmacro{\pigMethodPlotHeight}{\pigMethodAxisHeight+45pt}%
  \pgfmathsetmacro{\pigSourceY}{-(\pigBarWidth+2*\pigBarGap)/\pigRowPitch}%
  \pgfmathsetmacro{\pigSourceYMin}{\pigSourceY-.64}%
  \pgfmathsetlengthmacro{\pigSourcePlotHeight}{(4.64-\pigSourceYMin)*\pigRowPitch+45pt}%
}

\newcommand{\tableinput}[1]{\csname @@input\endcsname #1 }

\newcommand{\Z}{\ensuremath{\mathbb{Z}}}

\usepackage{xspace}

\usepackage{accents}
\newcommand{\ubar}[1]{\underaccent{\bar}{#1}}
\newcommand{\aC}{\ensuremath{\mathcal{C}}}   
\newcommand{\C}{\ensuremath{C}}              
\newcommand{\TAL}{\ensuremath{Z}}            

\newcommand{\Demand}[1]{\ensuremath{g_{#1}}}
\newcommand{\StatD}[1]{\ensuremath{D_{#1}}}
\newcommand{\DemandMin}[1]{\ensuremath{\ubar{G}_{#1}}}
\newcommand{\DemandMax}[1]{\ensuremath{\bar{G}_{#1}}}
\newcommand{\Pot}[1]{\ensuremath{\theta_{#1}}}
\newcommand{\PotC}[1]{\ensuremath{\Pot{#1}^{\C}}}
\newcommand{\z}[2]{\ensuremath{z_{({#1},{#2})}}}

\newcommand{\Flow}[2]{\ensuremath{x_{({#1},{#2})}}}
\newcommand{\FlowC}[2]{\ensuremath{\Flow{#1}{#2}^{\C}}}
\newcommand{\FlowOverC}[3]{\ensuremath{o_{({#1},{#2})}^{#3}}}

\newcommand{\Tract}{\ensuremath{\textsc{Tract}}}
\newcommand{\Opt}{\ensuremath{\textsc{Opt}}}

\newcommand{\problemstatement}[4]{%
  \par\medskip
  \begingroup
  \setlength{\fboxsep}{0pt}%
  \noindent\colorbox{black!4}{%
    \begin{minipage}{\linewidth}
      \small
      \colorbox{black!14}{%
        \parbox{\linewidth}{\vspace{0.45em}\hspace{0.75em}\textsc{#1}\par\vspace{0.4em}}%
      }%
      \vspace{0.4em}
      \hspace*{0.75em}%
      \begin{minipage}{\dimexpr\linewidth-1.5em\relax}
        \renewcommand{\arraystretch}{1.18}%
        \begin{tabularx}{\linewidth}{@{}>{\bfseries}l@{\hspace{0.8em}}>{\raggedright\arraybackslash}X@{}}
          Instance. & #2 \\
          Feasible solution. & #3 \\
          Objective. & #4
        \end{tabularx}
      \end{minipage}%
      \vspace{0.4em}
    \end{minipage}%
  }%
  \endgroup
  \par\medskip
}

\AddToHook{env/table/begin}{\let\scriptsize\footnotesize}

\renewcommand{\FloatBarrier}{}
\usepackage{authblk}
\usetikzlibrary{svg.path}
\definecolor{orcidlogocol}{HTML}{A6CE39}
\tikzset{
    orcidlogo/.pic={
        \fill[orcidlogocol] svg{M256,128c0,70.7-57.3,128-128,128C57.3,256,0,198.7,0,128C0,57.3,57.3,0,128,0C198.7,0,256,57.3,256,128z};
        \fill[white] svg{M86.3,186.2H70.9V79.1h15.4v48.4V186.2z}
        svg{M108.9,79.1h41.6c39.6,0,57,28.3,57,53.6c0,27.5-21.5,53.6-56.8,53.6h-41.8V79.1z M124.3,172.4h24.5c34.9,0,42.9-26.5,42.9-39.7c0-21.5-13.7-39.7-43.7-39.7h-23.7V172.4z}
        svg{M88.7,56.8c0,5.5-4.5,10.1-10.1,10.1c-5.6,0-10.1-4.6-10.1-10.1c0-5.6,4.5-10.1,10.1-10.1C84.2,46.7,88.7,51.3,88.7,56.8z};
    }
}

\newcommand{\paperorcid}[1]{\begingroup\hypersetup{hidelinks}\href{https://orcid.org/#1}{%
  \begin{tikzpicture}[xscale=0.025,yscale=-0.025,transform shape]
    \pic{orcidlogo};
  \end{tikzpicture}}\endgroup}

\newtheorem{lemma}{Lemma}

\newtheorem{proposition}{Proposition}

\newtheorem{corollary}{Corollary}
\theoremstyle{remark}

\theoremstyle{plain}

\newcommand{\paperfrontmatter}[1]{%
  \hypersetup{pdftitle={\paperTitle},pdfauthor={\paperAuthor}}%
  \title{\paperTitle}%
  \author{\paperAuthor\,\protect\paperorcid{\paperORCID}}%
  \affil{\paperAffiliation\\
    \protect\href{mailto:donkiewicz@or.rwth-aachen.de}{donkiewicz@or.rwth-aachen.de}}%
  \date{}%
  \maketitle
  \begin{abstract}#1\end{abstract}
  \noindent\textbf{Keywords:} \paperKeywords
  \par\medskip
}
\newcommand{\paperappendix}{\clearpage}

\newcommand{\paperTitle}{Composition of Near-Planar Infrastructure Network Instances
with Certified Skewness Bounds}
\newcommand{\paperAuthor}{Tim Donkiewicz}
\newcommand{\paperAffiliation}{Chair of Operations Research, RWTH Aachen University\\
Kackertstra\ss e 7, 52072 Aachen, Germany}
\newcommand{\paperORCID}{0000-0002-5721-3563}
\newcommand{\paperKeywords}{network instance generation; graph composition;
skewness; maximum planar subgraph; infrastructure networks; computational
experiments; network optimization.}

\begin{document}
\paperfrontmatter{
Computational studies of optimization problems on infrastructure networks may require more suitable instances than are available.
We propose a method to obtain additional instances by extracting and joining connected parts of existing networks, reusing their topology and parameters.
Many infrastructure networks are near-planar, so added connections should preserve this property.
We measure near-planarity by skewness, the minimum number of edges whose removal makes a graph planar.
We derive bounds on the additional skewness caused by joining and give sufficient conditions under which none is added, and propose the corresponding strategy.
We evaluate our approach on 100 sets of network parts extracted from each of PGLib-OPF~\cite{pglib} and SNDlib~\cite{orlowski2010sndlib}.
Within time limits, we show that no skewness is added for about half of the resulting power transmission networks and all considered telecommunication networks.
Our method gives smaller mean upper bounds on added skewness than baselines that select vertices randomly, in order of decreasing degree, or using a Kamada--Kawai layout.
Extraction and composition may also make the resulting optimization instances infeasible.
To reduce this risk, we solve a simpler, tractable variant of the intended optimization problem during generation and restart if it is infeasible.
On the generated instances, we solve security-constrained optimal transmission switching and survivable network design problems in full, demonstrating that the instances vary in computational difficulty with runtimes ranging from seconds to hours.
The method thus extends the available benchmark sets with instances that inherit source data, have known skewness bounds, and cover a range of computational difficulty.
}

\section{Introduction}\label{sec:introduction}

Meaningful computational studies require a sufficiently large and diverse test set to support their conclusions~\cite{beiranvand2017best,demsar2006statistical,campelo2019sample}.
The number of instances provided by public libraries may be insufficient to reliably establish statistically significant performance differences between algorithms, leaving researchers to draw unreliable conclusions, or to obtain additional instances.
Additional instances can be obtained by generating them from scratch or reusing existing ones.
The difficulty in the former is choosing realistic parameters, which requires considerable domain knowledge.
The advantage in the latter is that the instances are already accepted by the community, and for certain types of problems, a combination of multiple instances may also be accepted if an appropriate combination strategy is applied.

\begin{figure}[tbp]
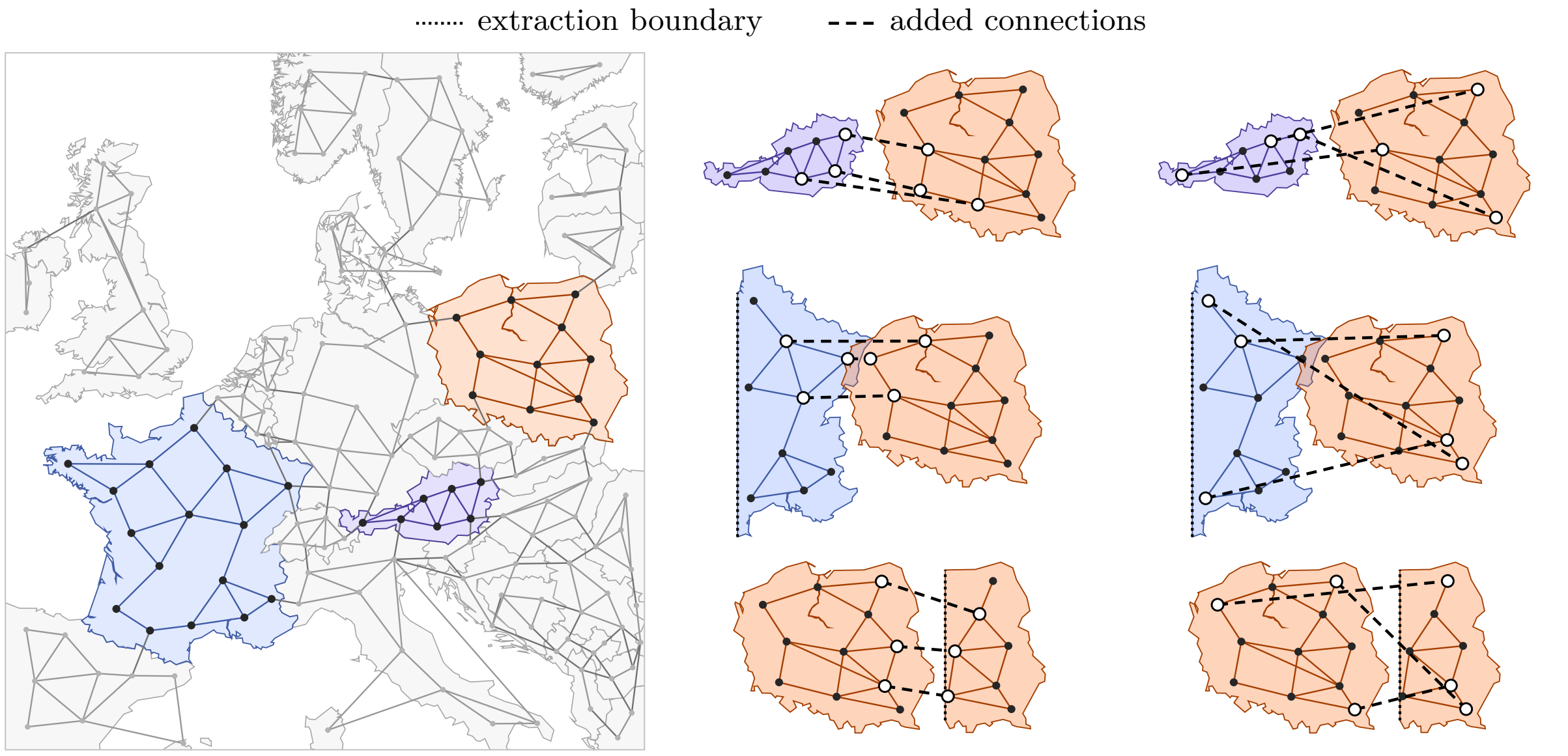

  \centering
  \input{resources/country-composition/colors.tex}
  \pgfmathsetlengthmacro{\countryfigureunit}{\linewidth/15.8}
  \tikzset{
    new connection/.style={draw=black,line width=0.85pt,dash pattern=on 3pt off 2pt},
    extraction boundary/.style={draw=black,line width=0.8pt,densely dotted}
  }
  {\footnotesize
    \tikz[baseline=-.5ex]{\draw[extraction boundary] (0,0)--(.5,0);}
    extraction boundary\qquad
    \tikz[baseline=-.5ex]{\draw[new connection] (0,0)--(.5,0);}
    added connections
  }
  \par\smallskip
  \begin{subfigure}[c]{0.42\linewidth}
    \centering
    \input{resources/country-composition/geography.tex}
    \caption{European power grid}
    \label{fig:input-networks-geography}
  \end{subfigure}\hfill
  \begin{subfigure}[c]{0.28\linewidth}
    \centering
    \input{resources/country-composition/good.tex}
    \caption{Planar compositions}
    \label{fig:input-networks-good}
  \end{subfigure}\hfill
  \begin{subfigure}[c]{0.28\linewidth}
    \centering
    \input{resources/country-composition/poor.tex}
    \caption{Nonplanar compositions}
    \label{fig:input-networks-bad}
  \end{subfigure}
  \caption{Schematic networks and connections. Country boundaries:
    \textcopyright{} EuroGeographics, via Eurostat/GISCO.}
  \label{fig:input-networks}
\end{figure}

We consider network-based optimization problems for real-world infrastructure networks, where a composition comprises merging the underlying networks, e.g., by adding an edge (see \Cref{fig:input-networks}).
Infrastructure networks are usually near-planar~\cite{barthelemy2011spatial,birchfield2021crossings,cotillasanchez2012comparing} (\Cref{fig:input-networks-geography}), so it is desirable that generated instances retain that property (\Cref{fig:input-networks-good}).
We quantify near-planarity by \emph{skewness}, the minimum number of edges that must be removed to make a graph planar~\cite[Def.~11]{liebers2001planarizing}.
A poor merge of two networks can increase skewness well beyond that already present in the extracted subgraphs (\Cref{fig:input-networks-bad}).
Avoiding this increase may require computing maximum planar subgraphs, which is computationally expensive.

We study how to compose network instances while controlling added skewness and reducing infeasibility rates.
Our contributions are the following.
First, we propose a generic composition framework.
It extracts connected subgraphs, retaining their node and edge parameters, and compensates flows removed at cut edges.
We join the subgraphs by adding coupling edges between selected vertices, called \emph{merge vertices}.
We solve a tractable auxiliary problem of the problem the instance is intended for, and only retain the instance if it is feasible.
Second, we derive bounds on the skewness added by composition and give conditions under which none is added.
We also bound changes in degree distribution, maximum degree, and average clustering.
Third, we design an exact vertex selection strategy and prove that it adds no skewness when maximum planar subgraphs are used.
We compare this strategy with alternatives that select vertices randomly, in order of decreasing degree, or based on a Kamada--Kawai layout~\cite{kamada1989algorithm}.
Fourth, we evaluate the construction on source instances stemming from PGLib-OPF~\cite{pglib} and SNDlib~\cite{orlowski2010sndlib}.
We then use the generated instances in experiments in security-constrained optimal transmission switching (SCOTS) and in survivable network design (SND).

Our proposed exact selection strategy gives smaller mean upper bounds on added skewness than the alternatives.
We show that no skewness is added for about half of the power transmission composites and all telecommunication composites.
The choice of merge vertices also affects properties beyond skewness, including degree distribution, graph diameter and connectivity, which we compare across the different strategies.
To examine whether the generated instances are useful for optimization experiments, we also evaluate their computational outcomes and solution times.
Runtimes range from seconds to hours in both applications, demonstrating differences in computational difficulty among the generated instances.

The proposed method therefore extends the available test sets with instances that retain characteristics of realistic infrastructure networks and span a range of computational difficulty.

\section{Related Work}\label{sec:related-work}

\paragraph*{Computational Evaluation: Best Practices.}
Guidelines for computational experiments emphasize representative, sufficiently large and diverse test sets~\cite{mcgeoch2012guide,beiranvand2017best,smithmiles2015generating}.
\citet{hooker1994needed,hooker1995testing} additionally suggests systematic variation of instance properties helps explain algorithm performance.
\citet{beiranvand2017best} recommend at least $20$ instances for a dependable conclusion and $100$ or more to detect performance differences between algorithm variants.
\citet{campelo2019sample} and \citet{campelo2020sample} provide methods to determine how many instances are needed to detect specified performance differences at a chosen statistical significance level.
For example, they calculate a requirement of $34$ instances for a paired $t$-test at a $5\%$ significance level~\cite[Section~5.1]{campelo2019sample}.
\citet{lodi2013performance} also recommend repeated runs on the same instance (e.g., with different seeds) to account for runtime variability in mixed-integer programming.

These recommendations motivate sufficiently large and diverse test sets that remain representative of the application.
For the infrastructure networks considered in this paper, this means retaining near-planarity while varying properties such as network size, connectivity and computational difficulty.

\paragraph*{Computational studies of security-constrained optimal power flow problems.}\label{sec:rw:modification}
Optimal power flow (OPF) problems minimize the cost of generation under the physical constraints of the power grid.
The security-constrained optimal power flow (SCOPF) problem extends the OPF by requiring feasible operation under contingency scenarios, e.g., single-edge failures.
PGLib-OPF provides a library of $66$ power system instances ranging from $3$ to $78{,}484$ buses~\cite{pglib}.
Due to that large range, the number of \emph{relevant} instances for a given problem formulation and algorithm (i.e., those that are neither trivial, nor too hard) can be much smaller.
Additionally, some instances do not admit a feasible solution to SCOPF or its variants, motivating researchers to modify them to obtain feasible instances.

\Cref{tab:modifications} summarizes our review of the Web of Science Core Collection~\cite{webofscience} query \texttt{TS=("optimal power flow") AND PY=(2003--2023)}, executed on 2024-04-09.
After sorting by citation count, we screened the $50$ most-cited records and retained $18$ computational studies that impose contingency constraints or analyze the computed solutions under contingencies.


\newsavebox{\bxfull}\newsavebox{\bxhalf}\newsavebox{\bxempty}
\newsavebox{\bxhalfq}
\newsavebox{\bxfullg}\newsavebox{\bxhalfg}\newsavebox{\bxemptyg}
\newsavebox{\bxhalfgq}
\tikzset{table symbol/.style={x=3mm,y=3mm,line width=0.4pt}}
\sbox{\bxfull}{\tikz[table symbol]{\filldraw[black] (0,0) rectangle (1,1);}}
\sbox{\bxhalf}{\tikz[table symbol]{\fill[black] (0,0) -- (1,1) -- (1,0) -- cycle;\draw[black] (0,0) rectangle (1,1);}}
\sbox{\bxhalfq}{\tikz[table symbol]{
    \fill[black] (0,0) -- (1,1) -- (1,0) -- cycle;
    \draw[black] (0,0) rectangle (1,1);
    \node[text=black,font=\sffamily\bfseries\fontsize{8}{8}\selectfont,inner sep=0pt,outer sep=0pt] at (0.5,0.5) {?};
    \begin{scope}
        \clip (0,0) -- (1,1) -- (1,0) -- cycle;
        \node[text=white,font=\sffamily\bfseries\fontsize{8}{8}\selectfont,inner sep=0pt,outer sep=0pt] at (0.5,0.5) {?};
    \end{scope}
}}
\sbox{\bxempty}{\tikz[table symbol]{\draw[black] (0,0) rectangle (1,1);}}
\sbox{\bxfullg}{\tikz[table symbol]{\filldraw[black] (0.5,0.5) circle[radius=0.5];}}
\sbox{\bxhalfg}{\tikz[table symbol]{\fill[black] (0.5,1) arc[start angle=90,end angle=-90,radius=0.5] -- cycle;\draw[black] (0.5,0.5) circle[radius=0.5];}}
\sbox{\bxemptyg}{\tikz[table symbol]{\draw[black] (0.5,0.5) circle[radius=0.5];}}
\sbox{\bxhalfgq}{\tikz[table symbol]{
    \fill[black] (0.5,1) arc[start angle=90,end angle=-90,radius=0.5] -- cycle;
    \draw[black] (0.5,0.5) circle[radius=0.5];
    \node[text=black,font=\sffamily\bfseries\fontsize{8}{8}\selectfont,inner sep=0pt,outer sep=0pt] at (0.5,0.5) {?};
    \begin{scope}
        \clip (0.5,0) rectangle (1,1);
        \node[text=white,font=\sffamily\bfseries\fontsize{8}{8}\selectfont,inner sep=0pt,outer sep=0pt] at (0.5,0.5) {?};
    \end{scope}
}}
\newcommand{\fullrect}{\usebox{\bxfull}}
\newcommand{\halfrect}{\usebox{\bxhalf}}
\newcommand{\halfrectq}{\usebox{\bxhalfq}}
\newcommand{\emptyrect}{\usebox{\bxempty}}
\newcommand{\fullcircleg}{\usebox{\bxfullg}}
\newcommand{\halfcircleg}{\usebox{\bxhalfg}}
\newcommand{\emptycircleg}{\usebox{\bxemptyg}}
\newcommand{\halfcirclegq}{\usebox{\bxhalfgq}}
\newcommand{\incompleteinfo}{\textsuperscript{\sffamily\bfseries ?}}

\newcommand{\namecite}[1]{\citet{#1}}

\begin{table*}[t]
    \centering
    \caption{Test case use in power flow and transmission switching studies. Base networks: \fullrect{} unmodified, \halfrect{} modified, \emptyrect{} unused. Numbers denote base-network bus counts. Enforced contingencies within each study's outage model: \fullcircleg{} all, \halfcircleg{} a subset, \emptycircleg{} none. A question mark indicates undisclosed or incomplete modification information.}
    \label{tab:modifications}
    \small
    \setlength{\tabcolsep}{1.5pt}
    \begin{tabular}{@{}l*{6}{c}lc@{}}
    \toprule
                                            & \multicolumn{7}{c}{\textbf{Base network}} & \multirow{2}{*}{\textbf{Outages}} \\
    \cmidrule(lr){2-8}
                                            & \textbf{14}& \textbf{30}& \textbf{57}& \textbf{73}&\textbf{118}&\textbf{300}& \textbf{Other}                          & \\
    \midrule
    \namecite{fisher2008optimal}            & \emptyrect & \emptyrect & \emptyrect & \emptyrect & \halfrect  & \emptyrect &                                         & \emptycircleg \\
    \namecite{hedman_optimal_2009}          & \emptyrect & \emptyrect & \emptyrect & \halfrect  & \halfrect  & \emptyrect &                                         & \halfcircleg  \\
    \namecite{khodaei2010transmission}      & \emptyrect & \emptyrect & \emptyrect & \emptyrect & \halfrect  & \emptyrect & Custom 6-bus                            & \halfcirclegq \\
    \namecite{o2010economic}                & \emptyrect & \emptyrect & \emptyrect & \halfrect  & \emptyrect & \emptyrect &                                         & \fullcircleg  \\
    \namecite{barrows2012computationally}   & \emptyrect & \emptyrect & \emptyrect & \halfrect  & \emptyrect & \emptyrect &                                         & \fullcircleg  \\
    \namecite{khanabadi2012optimal}         & \emptyrect & \emptyrect & \halfrectq & \emptyrect & \emptyrect & \fullrect  &                                         & \halfcircleg  \\
    \namecite{khodaei2012security}          & \emptyrect & \emptyrect & \emptyrect & \emptyrect & \halfrect  & \emptyrect &                                         & \halfcircleg  \\
    \namecite{ruiz2012reduced}              & \emptyrect & \emptyrect & \emptyrect & \emptyrect & \halfrectq & \emptyrect &                                         & \halfcirclegq \\
    \namecite{heidarifar2014optimal}        & \halfrect  & \emptyrect & \emptyrect & \emptyrect & \emptyrect & \emptyrect &                                         & \halfcircleg  \\
    \namecite{chiang2015solving}            & \emptyrect & \halfrectq & \halfrectq & \emptyrect & \halfrect  & \halfrect  & 3-bus; Iceland\incompleteinfo; grids\incompleteinfo & \halfcirclegq \\
    \namecite{ding2016robust}               & \emptyrect & \halfrect  & \emptyrect & \emptyrect & \halfrect  & \emptyrect &                                         & \fullcircleg  \\
    \namecite{khanabadi2018decentralized}   & \emptyrect & \emptyrect & \emptyrect & \emptyrect & \halfrectq & \emptyrect & 6-bus; 472-bus\incompleteinfo            & \halfcirclegq \\
    \namecite{behnia2019integrated}         & \emptyrect & \emptyrect & \emptyrect & \emptyrect & \emptyrect & \emptyrect & Modified IEEE-24                        & \fullcircleg  \\
    \namecite{saavedra_day-ahead_2020}      & \emptyrect & \emptyrect & \emptyrect & \emptyrect & \halfrect  & \halfrect  & Custom 4-bus                            & \halfcircleg  \\
    \namecite{han2021ensuring}              & \emptyrect & \fullrect  & \emptyrect & \emptyrect & \emptyrect & \emptyrect & Germany\incompleteinfo                  & \halfcircleg  \\
    \namecite{heidarifar2021optimal}        & \halfrect  & \emptyrect & \halfrect  & \emptyrect & \halfrect  & \emptyrect &                                         & \halfcircleg  \\
    \namecite{li2021connectivity}           & \emptyrect & \emptyrect & \emptyrect & \emptyrect & \emptyrect & \fullrect  & Modified PJM-5; French 1888-bus          & \halfcirclegq \\
    \namecite{han2022formulating}           & \fullrect  & \fullrect  & \fullrect  & \emptyrect & \emptyrect & \emptyrect & Germany (50Hertz)                       & \halfcircleg  \\
    \bottomrule
    \end{tabular}
\end{table*}

Most reviewed publications use only one to three base networks.
Some studies leave the exact contingency sets or details of other instance modifications undisclosed, limiting reproducibility and comparability.

The use of so few base networks and extensive, sometimes undisclosed modifications highlights the community's need for additional instances at relevant levels of computational difficulty to support meaningful computational evaluations.

\paragraph*{Computational studies of survivable network design problems.}\label{sec:rw:snd}
Network design problems minimize the cost of installing edge capacities to accommodate commodity flows between their sources and sinks.
A possible survivability requirement is that flows can be rerouted (arbitrarily) to fulfill all demands after any single-edge failure.
Available libraries provide topologies of real telecommunication networks~\cite{knight2011topologyzoo,spring2002rocketfuel} and graphs from a wide range of other applications~\cite{leskovec2014snap,rossi2015networkrepository}.
These libraries supply only the graph, requiring researchers to synthesize the parameters of the survivable network design problem themselves.
One of the most widely used libraries that does include parameters is SNDlib~\cite{orlowski2010sndlib}, whose original release contains $22$ network topologies.

\paragraph*{Planarization and maximum planar subgraphs.}
Computing a maximum planar subgraph, a planar subgraph with the largest possible number of edges, is NP-hard~\cite{yannakakis1978nodeedge}.
Exact approaches include integer programming formulations solved by branch-and-cut~\cite{junger1996maximum,chimani2018cycles,chimani2018exact}.
\citet{liebers2001planarizing} surveys the broader planarization literature.

\subsection*{Instance Generation Strategies}\label{sec:rw:strategies}
We categorize approaches to generating network instances into five groups.

\emph{Generation} builds an instance from scratch by sampling from a model.
It requires significant domain knowledge to design a generator that produces realistic instances.
Application-specific generators produce water distribution~\cite{decorte2014hydrogen}, heating~\cite{buchholz2026berlin}, and transportation networks~\cite{bonnet2025ssndp,kovaleva2026logistic}.
Other approaches construct power system models from geographic and infrastructure data~\cite{birchfield2017synthetic,horsch2018pypsa,medjroubi2017opendata}.
Graph generators sample from random graph models to produce large networks with specified structural properties, such as degree distributions and clustering coefficients~\cite{hamann2016emlfr,blasius2019girg,funke2017kagen}.

\emph{Perturbation} derives new instances by changing the parameters or topology of an existing instance.
Examples include varying loads, power generation, and topology in power systems~\cite{lovett2024opfdata,rivera2025pfdelta}, and water demands, pipe diameters, and pump speeds in water networks~\cite{truong2025ditec}.
OPF-Learn~\cite{joswigjones2022opflearn} generates feasible load profiles for the alternating-current (AC) optimal power flow problem for a given network.
Multiscale Planar Graph Generation~\cite{chauhan2019multiscale} modifies nodes and edges of a given planar graph to generate structurally similar planar graphs.
ReCoN~\cite{staudt2017recon} and Global Curveball~\cite{carstens2018curveball} randomly replace edges while preserving each node's degree.

\emph{Learning} fits a generative model to existing graphs and samples new graphs from it~\cite{leskovec2010kronecker,you2018graphrnn,bojchevski2018netgan}.
For power systems, \citet{xiao2026operable} learn from instances with feasible AC power flows to generate new feasible instances.

\emph{Slicing} extracts smaller subinstances from one large instance.
\citet[Section~IV.A]{li2026sampfa} extract connected subgraphs of power networks and replace flows on removed connections with equivalent loads.

\emph{Composition} joins several existing instances into a new, larger one.
\citet{grigg1999rts96} construct the IEEE 73-bus system by connecting three copies of the IEEE 24-bus system.
They add one bus connected to an existing bus by a phase-shifting transformer.
\citet{donkiewicz2026adaptive} extracts subgraphs from SNDlib instances and connects high-degree nodes to generate further SND instances.
Classical work in graph theory studies how merging graphs by identifying selected vertices affects the structural properties of the resulting graph~\cite{wagner1937clique,decker1981genus}.

\subsection{Positioning of this work}\label{sec:rw:positioning}
Our method reuses the topology and parameters of existing infrastructure network instances, reducing the need to synthesize application data.
We select source instances within an application domain from existing, accepted libraries.
From those source instances, we extract subinstances and join them with coupling edges, retaining internal edges.
This allows one composite to combine characteristics of different source networks.
Multiscale Planar Graph Generation and ReCoN generate variants of a single source graph without assigning application parameters~\cite{chauhan2019multiscale,staudt2017recon}.
Relying on one source can bias a test set toward that network.

In graph composition, arbitrary choices of where to put coupling edges can produce unrealistic infrastructure networks by disrupting their characteristic near-planarity (\Cref{fig:input-networks-bad}).
We extend and compare to earlier composition approaches~\cite{grigg1999rts96,donkiewicz2026adaptive} by using planar subgraphs to select merge vertices.
We derive bounds on added skewness and prove that this strategy adds none when the planar subgraphs are maximum.
We also bound changes in degree distribution, maximum degree, and average clustering.

Like \citet{li2026sampfa}, we compensate removed flows, but also allow the compensation to vary to improve feasibility after composition.
We establish sufficient conditions for preserving feasibility of an application-specific auxiliary problem during extraction and composition (\Cref{prop:feasibility,cor:potential-join}).
We solve this auxiliary problem for each composite.
We evaluate the generated instances' structural properties and computational difficulty for SCOTS and SND.

\section{Problem Setting and Notation}\label{sec:problem-setting}
To establish a generic framework for the composition of instances for network-based optimization problems, we first define an abstract network instance and the associated optimization problem (\Cref{sec:problem:abstract,sec:problem:opt}).
We define graph-theoretic properties used to describe the graphs underlying a network instance and how they change when composing multiple graphs (\Cref{sec:problem:structure,sec:problem:composition,sec:problem:map}).
Finally, we introduce the two concrete problems used for the computational evaluation: security-constrained optimal transmission switching and its direct-current (DC) optimal power flow feasibility check (\Cref{sec:problem:scots}), and min-cost multicommodity flow as the feasibility check for SND (\Cref{sec:problem:snd}).

Note that throughout this paper, we will be referring to \emph{(sub-)instances} whenever talking about the graph together with node and edge parameters and possibly some metainformation required for the optimization problem.
We refer to \emph{(sub-)graphs} whenever the described properties or algorithms do not access any of the instance data, but only require the graph.

\subsection{Abstract network instance}\label{sec:problem:abstract}

A \emph{network instance} is a tuple $\mathcal{N} = (V, E, K, \mathcal{S}, \mathcal{P})$.
We write $G(\mathcal{N})=(V,E)$ for its underlying graph.
The pair $(V,E)$ is a graph whose edges are written in a fixed reference direction.
For a directed model, this is the direction of the edge.
For an undirected model, it is only a notational orientation, and the application determines whether flow opposite to it is allowed.
The set $K$ contains the \emph{commodities} and may be empty or contain only one element.
A commodity is a class of flow with its own supply and demand.
The remaining components are families indexed by $K$.
$\mathcal{S} = \{S_k\}_{k\in K}$ collects all supply and demand values attached to a node.
For each $k\in K$, $S_k:V\to\mathbb{R}^2$ is a function whose value $S_{k1}(v)$ is the lower bound and $S_{k2}(v)$ is the upper bound for supply or demand (depending on whether it assumes a positive or negative value) at node $v$.
In particular, when $S_{k1}(v)=S_{k2}(v)$, the supply or demand at $v$ is fixed to that value. Otherwise, it is a variable within these bounds.
If $f_k$ and $s_k$ denote realized flows and supplies where $s_k(v)\in[S_{k1}(v), S_{k2}(v)]$, we impose flow conservation
\[
  \sum_{(v,w)\in E}f_k(v,w)
  =\sum_{(u,v)\in E}f_k(u,v)+s_k(v),
\]
i.e., total outgoing flow equals total incoming flow plus the supply attached to the node.
$\mathcal{P} = \{P_k\}_{k\in K}$ consists of node- and edge-attached parameter data per commodity.
This includes, e.g., capacities, physical attributes (susceptances and angle limits in power systems; routing information in SND) and objective coefficients (generation cost functions in power systems; unit routing costs in SND).
When $|K| = 1$, the commodity index is dropped from the notation.
Several commodities may share a single attribute; e.g., in the considered SND all commodities share the edge capacity (\Cref{sec:problem:snd}).

We call a finite set of network instances that the community accepts as realistic a \emph{library}, denoted $\mathfrak{N} = \{\mathcal{N}_1, \mathcal{N}_2, \dots\}$.
The composition pipeline of \Cref{sec:pipeline} draws its source instances from such a library.
When source instances are composed, an edge added between two different source graphs is a \emph{coupling edge}, and a vertex incident with such an edge is a \emph{merge vertex}.

\subsection{Abstract optimization problem and feasibility check}\label{sec:problem:opt}

The same network instance may be used for different variants of an optimization problem.
Here, $\Opt(\mathcal{N})$ denotes the fixed optimization problem for which we generate instances.
$\Tract(\mathcal{N})$ denotes an application-specific, tractable auxiliary optimization problem used during instance generation.
It serves as a proxy for the suitability of $\mathcal{N}$ for $\Opt(\mathcal{N})$ and need not certify feasibility or infeasibility for $\Opt(\mathcal{N})$.
We solve $\Tract(\mathcal{N})$ for each selected source instance and composed instance (see \Cref{sec:pipeline}).

\subsection{Skewness and other structural properties of graphs}\label{sec:problem:structure}

To describe the structural effects of composition, we define the graph-theoretic properties for which we give identities, bounds or empirical measurements in the evaluation.
For these structural properties, we consider the simple undirected graph $G$ obtained by ignoring edge orientations, removing self-loops and replacing parallel edges by a single edge.
Further, let $N_G(v)=\{u : uv\in E(G)\}$ be the neighborhood of $v$ and $d_G(v)=|N_G(v)|$ its degree.

Our primary graph property is the \emph{skewness}
\[
  \begin{aligned}
    \mathrm{sk}(G)&=\min\{|F| : F\subseteq E(G),\ G-F\text{ is planar}\}\\
    &=|E(G)|-\mathrm{mps}(G),
  \end{aligned}
\]
where $\mathrm{mps}(G)=\max\{|E'| : E'\subseteq E(G),\ (V(G),E')\text{ is planar}\}$.
A subgraph attaining this maximum is called a \emph{maximum planar subgraph}; computing one is also called \emph{planarization}.
Throughout, we denote planar subgraphs by $H$.
Skewness counts exactly how many edges must be removed to recover planarity (see \Cref{fig:noncrossing-matching-joined} for an example), and is thus a measure of the graph's near-planarity.
Computing a maximum planar subgraph is NP-hard~\cite{yannakakis1978nodeedge}.
If a maximum planar subgraph computation is stopped before proving optimality, any planar subgraph found gives an upper bound on skewness.

Some graph characteristics remain unchanged except in the vicinity of the merge vertices.
Those can be compared between the source graphs and the composite graphs.
We record the \emph{degree distribution}
\[
  p_G(q)=\frac{|\{v\in V(G) : d_G(v)=q\}|}{|V(G)|}
\]
and the maximum degree $\Delta(G)=\max_v d_G(v)$.

The \emph{local clustering coefficient} of $v$ is
\[
  c_G(v)=
  \begin{cases}
    \displaystyle\frac{2|\{uw\in E(G) : u,w\in N_G(v)\}|}{d_G(v)(d_G(v)-1)}, & d_G(v)\geq 2,\\[4pt]
    0, & d_G(v)<2
  \end{cases}.
\]
In other words, $c_G(v)$ is the fraction of possible triangles that include $v$.
We report the average clustering coefficient $\bar c(G)=\frac{1}{|V(G)|} \sum_v c_G(v)$.



We additionally consider characteristics whose changes cannot easily be bounded, because they depend on the choice of the merge vertices when composing two graphs.
Such characteristics can be compared against the full source library rather than the source graphs for the composite.
The \emph{diameter} of a connected graph gives the maximum distance between any two vertices:
\[
  \operatorname{diam}(G)=\max_{u,v\in V(G)}\operatorname{dist}_G(u,v),
\]
where $\operatorname{dist}_G(u,v)$ is the length of a shortest path between $u$ and $v$ in $G$.

Let $n=|V(G)|$ and let $0=\lambda_1(G)\leq\lambda_2(G)\leq\cdots\leq\lambda_n(G)$ be the eigenvalues of the Laplacian $L=D-A$, where $A$ is the adjacency matrix and $D$ is the diagonal matrix of vertex degrees, in the same vertex order as $A$.
The second eigenvalue $\lambda_2(G)$ is the \emph{algebraic connectivity}; it is positive exactly when $G$ is connected, and a small positive value indicates a bottleneck.

Finally, the \emph{degree assortativity} $r(G)$~\cite{newman2002assortative} is the correlation between endpoint degrees in the underlying simple undirected graph.
It is undefined when the variance of the endpoint degrees is zero.
When defined, it lies in $[-1,1]$: positive values indicate edges between vertices of similar degree, negative values indicate edges between vertices of dissimilar degree, and zero indicates no linear degree correlation.


\subsection{Skewness of composite graphs}\label{sec:problem:composition}
Consider connected source graphs $G_1,\dots,G_p$ and their \emph{disjoint union} $\mathbin{\dot\bigcup}_{i=1}^p G_i$, the graph whose vertex set and edge set are the disjoint unions of the $V(G_i)$ and the $E(G_i)$, so that the sources share no vertex or edge.
Adding a set $M$ of coupling edges, each between a vertex of one source and a vertex of another, yields the composite $C=(\mathbin{\dot\bigcup}_{i=1}^p G_i)+M$, where $+\,M$ denotes the addition of the edge set $M$.

\begin{figure}[tbp]
  \centering
  \begingroup
  \newcommand{\MatchingPanel}[1]{%
    \begin{tikzpicture}[x=.95cm,y=.95cm,
        every node/.style={font=\footnotesize},
        vertex/.style={circle,draw=black,fill=white,inner sep=0pt,
                       minimum size=4.5pt,line width=.5pt},
        boundary/.style={line width=.75pt},
        source edge/.style={line width=.45pt},
        omitted edge/.style={draw=pigSrc1,line width=.95pt,densely dotted,
                             preaction={draw=white,line width=3pt}},
        matching edge/.style={draw=pigMerge,line width=.85pt,
                              dash pattern=on 3pt off 2pt}]
      \ifnum#1=0
        \tikzset{x=1.4cm,y=1.2cm}
        \path[use as bounding box] (-1.8,-1.6) rectangle (5.4,2.15);
        \coordinate (a1) at (.1,1.2);
        \coordinate (a2) at (1,0);
        \coordinate (a3) at (.1,-1.2);
        \coordinate (a4) at (-.6,0);
        \coordinate (a5) at (-1.4,0);
        \coordinate (a6) at (.35,0);
        \coordinate (b1) at (3.5,1.2);
        \coordinate (b2) at (2.6,0);
        \coordinate (b3) at (3.5,-1.2);
        \coordinate (b4) at (4.2,0);
        \coordinate (b5) at (5,0);
        \coordinate (b6) at (3.25,0);
        \path[fill=black!7] (a6)--(a2)--(a3)--(a4)--cycle;
        \path[fill=black!7] (b1)--(b6)--(b2)--(b3)--(b4)--cycle;
        \node at (.15,-.55) {$f_1\ (4)$};
        \node at (3.45,-.55) {$f_2\ (5)$};
        \node at (-.55,.45) {$(3)$};
        \node at (-.65,-.4) {$(3)$};
        \node at (-.03,.4) {$(3)$};
        \node at (.52,.37) {$(3)$};
        \node at (4.15,.45) {$(3)$};
        \node at (4.25,-.4) {$(3)$};
        \node at (3.08,.37) {$(3)$};
        \node at (.75,1.55) {$(4)$};
        \node at (2.85,1.55) {$(4)$};
        \node at (.1,1.95) {$H_1=G_1-\{a_5a_6\}$};
        \node at (3.5,1.95) {$H_2=G_2$};
      \else
        \path[use as bounding box] (-1.8,-2.3) rectangle (5.4,2.4);
        \coordinate (a1) at (-.2,.55);
        \coordinate (a2) at (1,.4);
        \coordinate (a3) at (1,-.4);
        \coordinate (a4) at (-1,-1.2);
        \coordinate (a5) at (0,-.25);
        \coordinate (a6) at (-1,1.2);
        \coordinate (b1) at (4.35,0);
        \coordinate (b2) at (2.6,.4);
        \coordinate (b3) at (2.6,-.4);
        \coordinate (b4) at (3.5,-1.2);
        \coordinate (b5) at (3.75,-.25);
        \coordinate (b6) at (3.5,1.2);
        \path[fill=black!7,even odd rule]
          (-1.45,-1.5) rectangle (4.85,1.5)
          (a6)--(a2)--(a3)--(a4)--cycle
          (b1)--(b6)--(b2)--(b3)--(b4)--cycle;
        \draw[draw=black!25,densely dotted,line width=.3pt]
          (-1.45,-1.5) rectangle (4.85,1.5);
        \draw[matching edge] (a6)--(b6) (a4)--(b4);
        \node at (0,2.15) {$\mathrm{sk}(G_1)=1$};
        \node at (3.5,2.15) {$\mathrm{sk}(G_2)=0$};
        \ifnum#1=1
          \draw[matching edge] (a2)--(b2) (a3)--(b3);
          \node at (1.8,.95) {$M$};
          \node at (1.8,-2.15) {$\mathrm{sk}(C)=\mathrm{sk}(G_1)+\mathrm{sk}(G_2)=1$};
        \else
          \draw[matching edge] (a2)--(b3);
          \draw[matching edge,preaction={draw=black!7,line width=3pt}]
            (a3)--(b2);
          \node at (1.8,.95) {$M'$};
          \node at (1.8,-2.15) {$\mathrm{sk}(C')=\mathrm{sk}(G_1)+\mathrm{sk}(G_2)+1=2$};
        \fi
      \fi
      \draw[boundary,draw=pigSrc1] (a6)--(a2)--(a3)--(a4)--cycle;
      \draw[boundary,draw=pigSrc2] (b1)--(b6)--(b2)--(b3)--(b4)--cycle;
      \draw[source edge,draw=pigSrc1]
        (a1)--(a2) (a1)--(a5)--(a3) (a4)--(a5) (a4)--(a1)--(a6);
      \draw[source edge,draw=pigSrc2]
        (b1)--(b2) (b1)--(b5)--(b3) (b4)--(b5);
      \ifnum#1>0
        \draw[omitted edge] (a5)--(a6);
      \fi
      \foreach \j in {2,3,4,6} {
        \node[vertex,draw=pigSrc1,fill=pigSrc1] at (a\j) {};
        \node[vertex,draw=pigSrc2,fill=pigSrc2] at (b\j) {};
      }
      \foreach \j in {1,5} {
        \node[vertex,draw=pigSrc1] at (a\j) {};
        \node[vertex,draw=pigSrc2] at (b\j) {};
      }
      \ifnum#1=0
        \node[above=3pt] at (a1) {$a_1$};
        \node[right=2pt] at (a2) {$a_2$};
        \node[below=3pt] at (a3) {$a_3$};
        \node[above left=1pt] at (a4) {$a_4$};
        \node[left=2pt] at (a5) {$a_5$};
        \node[below=2pt] at (a6) {$a_6$};
        \node[above=3pt] at (b1) {$b_1$};
        \node[left=2pt] at (b2) {$b_2$};
        \node[below=3pt] at (b3) {$b_3$};
        \node[above right=1pt] at (b4) {$b_4$};
        \node[right=2pt] at (b5) {$b_5$};
        \node[below=2pt] at (b6) {$b_6$};
      \else
        \node[below right=3pt] at (a1) {$a_1$};
        \node[above right=2pt] at (a2) {$a_2$};
        \node[below right=2pt] at (a3) {$a_3$};
        \node[below=3pt] at (a4) {$a_4$};
        \node[below right=1pt] at (a5) {$a_5$};
        \node[above=3pt] at (a6) {$a_6$};
        \node[right=2pt] at (b1) {$b_1$};
        \node[above left=2pt] at (b2) {$b_2$};
        \node[below left=2pt] at (b3) {$b_3$};
        \node[below=3pt] at (b4) {$b_4$};
        \node[below left=3pt] at (b5) {$b_5$};
        \node[above=3pt] at (b6) {$b_6$};
      \fi
    \end{tikzpicture}%
  }
  \begin{subfigure}[t]{.85\linewidth}
    \centering
    \MatchingPanel{0}
    \caption{Planar embeddings. Faces with most vertices in grey.}
    \label{fig:noncrossing-matching-faces}
  \end{subfigure}

  \medskip
  \begin{subfigure}[t]{.49\linewidth}
    \centering
    \MatchingPanel{1}
    \caption{Noncrossing matching.}
    \label{fig:noncrossing-matching-joined}
  \end{subfigure}\hfill
  \begin{subfigure}[t]{.49\linewidth}
    \centering
    \MatchingPanel{2}
    \caption{Crossing matching.}
    \label{fig:noncrossing-matching-crossed}
  \end{subfigure}
  \endgroup
  \caption{$H_1,H_2$ are maximum planar subgraphs. Numbers in parentheses count distinct boundary vertices. Shading marks the chosen faces before joining. The dotted edge belongs to $G_1$. Dashed edges match the filled vertices.}
  \label{fig:noncrossing-matching}
\end{figure}
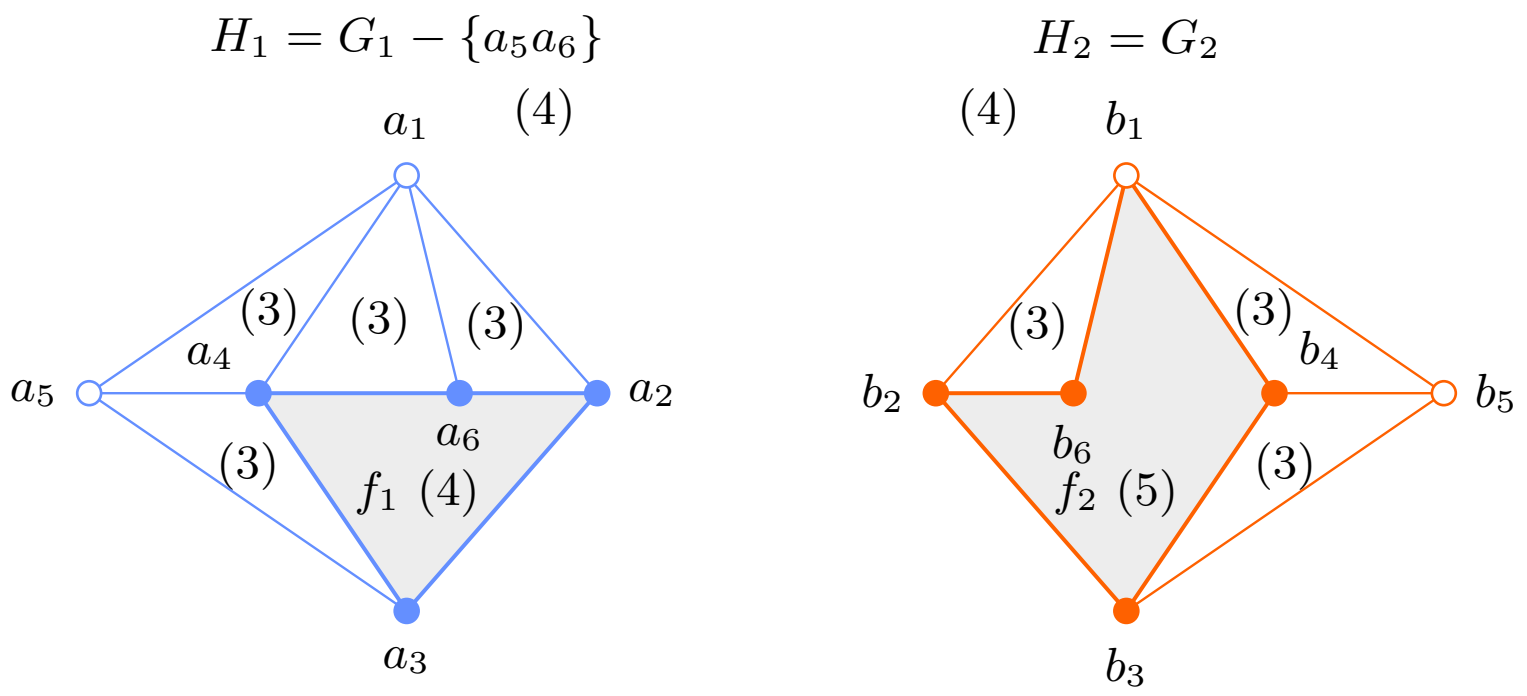
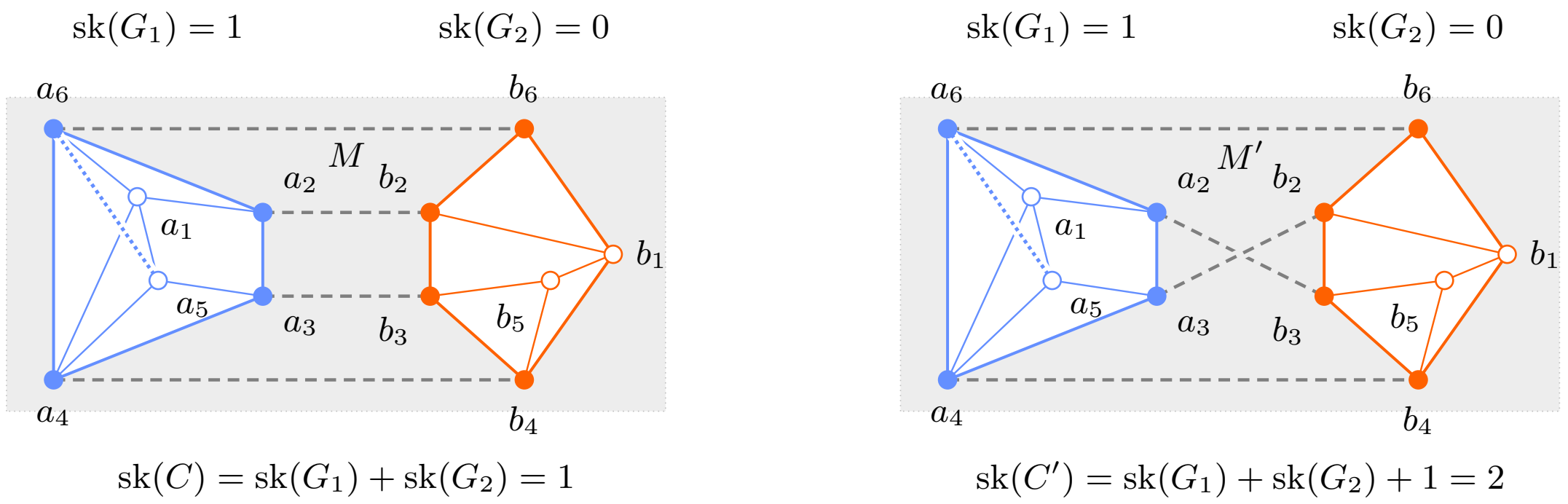

A \emph{planar embedding} of a graph is a drawing in the plane in which no two edges cross.
The embedded vertices and edges divide the plane into connected regions called \emph{faces} (see \Cref{fig:noncrossing-matching-faces}).
Exactly one face is unbounded and is called the \emph{outer face}.
An equivalent planar embedding can make any chosen face the outer face (see \Cref{fig:noncrossing-matching-faces} vs. \Cref{fig:noncrossing-matching-joined}, where $f_1$ and $f_2$ have been moved to the outside).
Following the boundary of a face $f$ gives a sequence of incident vertices and edges that returns to its starting vertex.
We denote this boundary walk by $W_f$, in which a vertex or edge may occur more than once (e.g., when a vertex of degree one lies on the boundary).
Let $V(f)$ be the set of distinct vertices occurring in $W_f$.
We order the vertices in $V(f)$ by their first occurrences in $W_f$.

A \emph{matching} between distinct vertices incident with chosen faces of two planar embeddings is called \emph{noncrossing} if, with the chosen faces placed as outer faces, the matching edges can be drawn without crossings (\Cref{fig:noncrossing-matching-joined}).
Not every vertex incident with either face must be used to give a valid matching.
A matching is \emph{crossing} if and only if it introduces skewness to the composite (\Cref{fig:noncrossing-matching-crossed}).

\begin{lemma}[Skewness bounds for composite graphs]\label{lem:skewness-bounds}
Suppose that the composite $C$ is connected.
Then the skewness of the composite $C=(\mathbin{\dot\bigcup}_{i=1}^p G_i)+M$ satisfies
\[
  \sum_{i=1}^p\mathrm{sk}(G_i)
  \leq \mathrm{sk}(C)
  \leq \sum_{i=1}^p\mathrm{sk}(G_i)+|M|-(p-1).
\]
\end{lemma}
\begin{proof}
The edges of $C$ split into the source edges and the coupling edges, so $|E(C)|=\sum_i|E(G_i)|+|M|$.
For the lower bound, let $H$ be a maximum planar subgraph of $C$.
Each restriction $H\cap G_i$ is planar and contains at most $\mathrm{mps}(G_i)$ edges.
Together with at most $|M|$ coupling edges, this gives $\mathrm{mps}(C)\leq\sum_i\mathrm{mps}(G_i)+|M|$ and hence $\mathrm{sk}(C)\geq\sum_i\mathrm{sk}(G_i)$.
For the upper bound, choose a maximum planar subgraph $H_i$ of each $G_i$.
Replace each source graph $G_i$ by a single vertex and retain the coupling edges between these vertices.
The resulting multigraph is connected and has $p$ vertices, so it has a spanning tree with $p-1$ edges.
Let $T\subseteq M$ be the corresponding set of coupling edges in $C$.
The graph $(\mathbin{\dot\bigcup}_i H_i)+T$ is planar: it joins vertex-disjoint planar graphs according to a tree, with every edge of $T$ acting as a bridge.
It omits $\sum_i\mathrm{sk}(G_i)$ source edges and $|M|-(p-1)$ coupling edges of $C$ (the ones not contained in the spanning tree) which proves the upper bound.
\end{proof}
For arbitrary (e.g., non-maximum) planar subgraphs that omit $d_i$ edges from $G_i$, the same argument gives $\mathrm{sk}(C)\leq\sum_i d_i+|M|-(p-1)$.

The following lemma characterizes when the lower bound is attained.
\begin{lemma}[Attainment of the lower skewness bound]\label{lem:floor}
The equality
\[
  \mathrm{sk}(C)=\sum_{i=1}^p\mathrm{sk}(G_i)
\]
holds if and only if there are maximum planar subgraphs $H_i\subseteq G_i$ for which $(\mathbin{\dot\bigcup}_{i=1}^p H_i)+M$ is planar.
In particular, for two sources, this condition holds whenever $M=\{a_jb_j : j=1,\dots,k\}$ is a noncrossing matching between distinct vertices incident with chosen faces of planar embeddings of maximum planar subgraphs $H_1$ and $H_2$.
Consequently, $k\leq\min\{|V(f_1)|,|V(f_2)|\}$ for the two chosen faces, i.e., $k$ can only be as large as the smaller of the two outer face sizes.
\end{lemma}
\begin{proof}
If $(\mathbin{\dot\bigcup}_iH_i)+M$ is planar, it is a planar subgraph of $C$ with $\sum_i\mathrm{mps}(G_i)+|M|$ edges.
Thus, $\mathrm{mps}(C)\geq\sum_i\mathrm{mps}(G_i)+|M|$.
The reverse inequality follows from the proof of \Cref{lem:skewness-bounds}, so equality holds.

Conversely, suppose the equality holds and let $H$ be a maximum planar subgraph of $C$.
By assumption, planarizing $C$ deletes exactly $\sum_i\mathrm{sk}(G_i)$ edges from its $\sum_i|E(G_i)|+|M|$ edges.
Therefore
\[
  \begin{aligned}
    |E(H)|&=\mathrm{mps}(C)\\
    &=\sum_i\bigl(|E(G_i)|-\mathrm{sk}(G_i)\bigr)+|M|\\
    &=\sum_i\mathrm{mps}(G_i)+|M|.
  \end{aligned}
\]
The right-hand side is the largest number of edges $H$ can have: at most $\mathrm{mps}(G_i)$ from each source $G_i$, plus the $|M|$ coupling edges.
Since $H$ has exactly that many edges, it contains all coupling edges and $\mathrm{mps}(G_i)$ edges from every source.
Thus, every $H_i=H\cap G_i$ is a maximum planar subgraph of $G_i$, and $H=(\mathbin{\dot\bigcup}_iH_i)+M$ is the required planar graph.
\end{proof}

\begin{lemma}[Boundary order matching]\label{lem:boundary-order}
Let $W_1$ and $W_2$ be outer face boundary walks of two vertex-disjoint planar graphs.
After fixing a starting vertex and traversal direction in each walk, enumerate their vertex occurrences as
\[
  W_i=(w^i_0,\dots,w^i_{\ell_i-1}),\qquad i\in\{1,2\},
\]
where $\ell_i$ is the number of vertex occurrences in one traversal, excluding the final repetition of the starting vertex.
Let $1\leq q\leq\min\{\ell_1,\ell_2\}$, and choose indices
\[
  0\leq \alpha_1<\cdots<\alpha_q<\ell_1
  \quad\text{and}\quad
  0\leq \beta_1<\cdots<\beta_q<\ell_2,
\]
i.e., some vertices may be skipped, but their order along each boundary walk must be preserved.
Suppose the vertices $a_j=w^1_{\alpha_j}$ are distinct and the vertices $b_j=w^2_{\beta_j}$ are distinct, i.e., at most one occurrence of each vertex is selected from each boundary walk.
Then the matching
\[
  M=\{a_jb_j : j=1,\dots,q\}
\]
is noncrossing.
\end{lemma}
\begin{proof}
We insert the matching edges $e_j=(a_j,b_j)$ in order.
An endpoint is called \emph{unused} until its matching edge is inserted; its position on the boundary walk is the one specified by the selected index, which may skip some intermediate nodes (e.g., because they are visited multiple times in the boundary walk).

Draw the two graphs disjointly with the chosen outer face nodes incident to the common outer face.
Observe that an edge can always be added without introducing crossings if its endpoints are incident to the same face of a planar embedding.
For $q\geq2$, we may need to mirror one drawing so that, after inserting $e_1$, the unused endpoints occur along the outer face boundary in the order $a_2,\dots,a_q,b_q,\dots,b_2$.
Thus, the first edge $e_1=a_1b_1$ can be drawn through the common outer face without crossings.
The induced face boundary walk becomes longer, traversing the new edge twice.
For $q=1$, the construction is complete after adding $e_1$.

Inductively, suppose that $e_1,\dots,e_j$ have been inserted without crossings, where $1\leq j<q$, and that the current graph has a face with a boundary walk that contains the unused endpoints in the order
\[
  a_{j+1},\dots,a_q,b_q,\dots,b_{j+1},
\]
possibly with intermediate vertices skipped according to the chosen indices, and the first edge $e_1$ and previously added edge $e_j$.
We add the next edge $e_{j+1}=a_{j+1}b_{j+1}$.
By construction, both $a_{j+1}$ and $b_{j+1}$ are incident with the same face of the current graph, and $e_{j+1}$ is noncrossing.
The new edge connects the first and last unused endpoints in the displayed order.
One resulting face is bounded by $e_j$ and $e_{j+1}$ with their endpoints, and other edges and vertices.
The other resulting face contains all remaining unused endpoints on its boundary in the same relative order.
Thus the induction hypothesis holds for $j+1$ whenever unused endpoints remain; after inserting $e_q$, the matching is complete without crossings.
\end{proof}

The same principle applies when more than two source graphs are composed sequentially.
For $r\geq2$, let $C_r$ denote the composite obtained after the $r-1$ joins of $G_1,\dots,G_r$ with coupling edge sets $M_2,\dots,M_r$.
\begin{corollary}[Sequential composition]\label{cor:sequential-floor}
Let $C_1=G_1$ and choose a maximum planar subgraph $H_1^C\subseteq C_1$. Then, compose with $G_2,\dots,G_p$ one at a time, with $C_r=(C_{r-1}\mathbin{\dot\cup}G_r)+M_r$.
At step $r$, let $H_r^G\subseteq G_r$ be a maximum planar subgraph and let $M_r$ be a noncrossing matching between distinct vertices incident with chosen faces of planar embeddings of $H_{r-1}^C$ and $H_r^G$.
If $H_r^C=H_{r-1}^C\cup H_r^G\cup M_r$ is retained for the next step, then the final composite $C_p$ satisfies
\[
  \mathrm{sk}(C_p)=\sum_{i=1}^p\mathrm{sk}(G_i).
\]
\end{corollary}
\begin{proof}
For the base case $r=1$, we have $C_1=G_1$, the chosen subgraph $H_1^C$ is maximum planar, and $\mathrm{sk}(C_1)=\mathrm{sk}(G_1)$.

For the induction step, let $2\leq r\leq p$ and assume that $H_{r-1}^C$ is a maximum planar subgraph of $C_{r-1}$ and that
\[
  \mathrm{sk}(C_{r-1})=\sum_{i=1}^{r-1}\mathrm{sk}(G_i).
\]
The noncrossing matching $M_r$ makes $H_r^C=H_{r-1}^C\cup H_r^G\cup M_r$ planar.
Applying \Cref{lem:floor} to $C_{r-1}$ and $G_r$ therefore gives
\[
  \mathrm{sk}(C_r)
  =\mathrm{sk}(C_{r-1})+\mathrm{sk}(G_r)
  =\sum_{i=1}^{r}\mathrm{sk}(G_i).
\]
The forward direction in the proof of \Cref{lem:floor} also shows that $H_r^C$ is a maximum planar subgraph of $C_r$.
Thus both induction claims hold for $r$, and setting $r=p$ proves the result.
\end{proof}

Writing $M=\bigcup_{r=2}^p M_r$ for the set of all coupling edges,
\Cref{lem:skewness-bounds} shows that the additional skewness
$\mathrm{sk}(C)-\sum_i\mathrm{sk}(G_i)$ of a sequential composition of
$p$ subgraphs lies between $0$ and $|M|-(p-1)$.
By \Cref{cor:sequential-floor}, it is zero (no skewness is added) if each join uses a noncrossing matching between vertices on chosen faces of maximum planar subgraphs.
\Cref{fig:noncrossing-matching-joined} illustrates such a join.

Suppose instead that only planar subgraphs $H_i'\subseteq G_i$, not necessarily maximum, are available.
Maintaining these planar subgraphs and using a noncrossing matching at every sequential composition step keeps their union with the coupling edges planar, by the same induction as in \Cref{cor:sequential-floor}.
The added skewness is bounded by the sum of the differences in edge counts between maximum planar subgraphs of $G_i$ and the chosen subgraphs $H_i'$.
\begin{lemma}[Skewness bound for non-maximum planar subgraphs]\label{lem:suboptimal-mps}
If $H_i'\subseteq G_i$ are planar subgraphs for which $(\mathbin{\dot\bigcup}_{i=1}^p H_i')+M$ is planar, then
\[
  \sum_{i=1}^p\mathrm{sk}(G_i)\;\le\;\mathrm{sk}(C)\;\le\;\sum_{i=1}^p\mathrm{sk}(G_i)+\sum_{i=1}^p\bigl(\mathrm{mps}(G_i)-|E(H_i')|\bigr).
\]
\end{lemma}
\begin{proof}
The lower bound is \Cref{lem:skewness-bounds}, which holds for any set of coupling edges.
For the upper bound, consider the combined graph consisting of the (non-maximum) planar subgraphs together with the coupling edges, $(\mathbin{\dot\bigcup}_i H_i')+M$.
It is planar by hypothesis and has $\sum_i|E(H_i')|+|M|$ edges.
This gives $\mathrm{mps}(C)\ge\sum_i|E(H_i')|+|M|$.
Since $|E(C)|=\sum_i|E(G_i)|+|M|$, the $|M|$ terms cancel out, giving
\begin{align*}
  \mathrm{sk}(C)&=|E(C)|-\mathrm{mps}(C)\\
  &\le\sum_i\bigl(|E(G_i)|-|E(H_i')|\bigr) \\
  &=\sum_i\bigl(\underbrace{|E(G_i)|-\mathrm{mps}(G_i)}_{\mathrm{sk}(G_i)}\bigr)\\
  &\quad+\sum_i\bigl(\underbrace{\mathrm{mps}(G_i)-|E(H_i')|}_{\text{difference in edge counts}}\bigr) \\
  &=\sum_i\mathrm{sk}(G_i)+\sum_i\bigl(\mathrm{mps}(G_i)-|E(H_i')|\bigr).
\end{align*}
\end{proof}
The additional skewness above the lower bound is thus at most the sum of the missing edge set cardinalities of the chosen subgraphs, $\sum_i\bigl(\mathrm{mps}(G_i)-|E(H_i')|\bigr)$, and this upper bound is zero when each $H_i'$ is a maximum planar subgraph, yielding \Cref{lem:floor}.

\subsection{Structural properties of composite graphs}\label{sec:problem:map}
We identify which structural properties of the original graphs remain unchanged and bound how much others can change under composition.

\begin{proposition}[Structural effects of composition]\label{prop:preserved}
  Let $C=(\mathbin{\dot\bigcup}_{i=1}^pG_i)+M$ be the composite graph.
  For $v\in V(C)$, let $d_M(v)$ be the number of coupling edges incident with $v$, let $B=\{v\in V(C) : d_M(v)>0\}$ be the set of merge vertices, and let $\Delta(M)=\max_v d_M(v)$.
  For (iv), additionally assume the sequential construction $C_1=G_1$ and $C_r=(C_{r-1}\mathbin{\dot\cup}G_r)+M_r$, where each $M_r$ is a matching between $V(C_{r-1})$ and $V(G_r)$.
  Then the following hold:
  \begin{enumerate}
    \item[(i)] $|V(C)|=\sum_i|V(G_i)|$ and $|E(C)|=\sum_i|E(G_i)|+|M|$.
    \item[(ii)] Let
    \[
      p_{\mathrm{src}}(q)
      =\frac{1}{|V(C)|}\sum_i |V(G_i)|\,p_{G_i}(q)
    \]
    be the degree distribution of the disjoint source graphs before composition. Then, for every $q$,
    \[
      |p_C(q)-p_{\mathrm{src}}(q)|
      \leq \frac{|B|}{|V(C)|},
    \]
    i.e., for any degree $q$, the fraction of degree-$q$ vertices changes by at most the fraction of merge vertices.
    \item[(iii)] The maximum degree satisfies
    \[
      \max\!\left\{\max_i\Delta(G_i),\Delta(M)\right\}
      \leq \Delta(C)
      \leq \max_i\Delta(G_i)+\Delta(M),
    \]
    i.e., composition cannot reduce the largest source maximum degree and can increase it by at most the maximum coupling degree.
    \item[(iv)]
    Any vertex $v\notin B$ satisfies $c_C(v)=c_{G_i}(v)$.
    Let
    \[
      \bar c_{\mathrm{src}}
      =\frac{1}{|V(C)|}\sum_i |V(G_i)|\,\bar c(G_i)
    \]
    be the average clustering coefficient over all source graphs before composition. Then
    \[
      0\leq \bar c_{\mathrm{src}}-\bar c(C)
      \leq \frac{|B|}{|V(C)|},
    \]
    i.e., the average clustering coefficient can decrease by at most the fraction of merge vertices.
  \end{enumerate}
\end{proposition}
\begin{proof}
The source graphs are vertex-disjoint and $M$ adds exactly $|M|$ edges, proving (i).

For $v\in V(G_i)$, counting its source and coupling edges gives
$d_C(v)=d_{G_i}(v)+d_M(v)$.
Only the vertices in $B$ change degree; hence, for each $q$, at most $|B|$ vertices switch their degree-$q$ class, proving the bound on $p_C(q)$ in (ii).
Taking maxima in the degree identity gives (iii).
If $v\notin B$, its neighborhood is unchanged, so $c_C(v)=c_{G_i}(v)$.

Suppose a triangle meets multiple source graphs, and let $r$ be the largest index for which it meets $G_r$.
Irrespective of whether it has one or two vertices in $G_r$, the two edges between $G_r$ and earlier source graphs are contained in $M_r$ and share an endpoint, contradicting that $M_r$ is a matching.
Thus no triangle in $C$ contains a coupling edge.
Therefore, composition leaves the number of triangles through each vertex unchanged while its degree can only increase, so $c_C(v)\leq c_{G_i}(v)$.
Since equality holds outside $B$ and all clustering coefficients lie in $[0,1]$, averaging proves (iv).


\end{proof}

\subsection{Maintaining feasibility through cut flow compensation}\label{sec:problem:feasibility}

When extracting subinstances from a network flow instance, the nodes incident to cut edges lose flow contributions from the removed edges.
This can make the subinstance infeasible, and make it more difficult for the composed instance to be feasible.

We apply cut flow compensation, adding the flow that previously crossed the cut to the supply at each node.
For ordinary network flow instances satisfying the hypotheses of \Cref{prop:feasibility}, cut flow compensation is sufficient to maintain feasibility for $\Tract$ on the subinstance. 
The final composed instance remains feasible for $\Tract$ as well if assigning zero flow to all coupling edges satisfies their bounds and all additional constraints.

Let $\mathcal{N}'$ be induced by $V'\subseteq V$, with $E'=E\cap(V'\times V')$, and let $\bar f_k$ be the flows from a feasible solution of $\Tract(\mathcal{N})$.
For $v\in V'$ and $k\in K$, the cut flow compensation is $v$'s net flow across the cut,
\begin{equation}\label{eq:cut-flow-compensation}
  \Delta s_k(v)
  =\sum_{\substack{(u,v)\in E\\u\notin V'}}\bar f_k(u,v)
  -\sum_{\substack{(v,w)\in E\\w\notin V'}}\bar f_k(v,w),
\end{equation}
which can be positive or negative.

\begin{proposition}[Cut flow compensation]\label{prop:feasibility}
Let $\bar{x}$ be a feasible solution for $\Tract(\mathcal{N})$, with realized flows $\bar f_k$ and supplies $\bar s_k$, and let $\Delta s_k$ be given by \eqref{eq:cut-flow-compensation}.
For every $v\in V'$ and $k\in K$, shift the supply limits in $\mathcal{N}'$ either to $S_{k1}(v)+\Delta s_k(v)$ and $S_{k2}(v)+\Delta s_k(v)$ (fixed compensation), or to $S_{k1}(v)+\min\{0,\Delta s_k(v)\}$ and $S_{k2}(v)+\max\{0,\Delta s_k(v)\}$ (variable compensation).
Suppose that every constraint of $\Tract(\mathcal{N}')$ other than flow conservation and supply limits is satisfied by the restriction of $\bar{x}$ together with the compensated supplies.
Then this restriction, with realized supplies $\bar s_k(v)+\Delta s_k(v)$, is feasible for $\Tract(\mathcal{N}')$.
\end{proposition}
\begin{proof}
For every $v\in V'$ and $k\in K$, feasibility of $\bar{x}$ in $\mathcal{N}$ gives
\[
  \sum_{(v,w)\in E}\bar f_k(v,w)
  =
  \sum_{(u,v)\in E}\bar f_k(u,v)+\bar s_k(v).
\]
Separating internal and cut edges and rearranging yields
\[
  \sum_{(v,w)\in E'}\bar f_k(v,w)
  =
  \sum_{(u,v)\in E'}\bar f_k(u,v)
  +\bar s_k(v)+\Delta s_k(v).
\]
Thus the restricted solution satisfies flow conservation with the compensated supplies.
Since $\bar s_k(v)\in[S_{k1}(v),S_{k2}(v)]$, the compensated supplies lie within the shifted limits for both compensations, and the hypothesis covers all remaining constraints.
\end{proof}

The above proposition guarantees feasibility of the subinstance under its hypotheses, but does not by itself guarantee feasibility of the composed instance, because additional constraints may be imposed on endpoints of coupling edges.
For example, in potential-based flow problems, each node $v$ has a continuous potential $\pi_v$, and the flow on an edge is determined by the difference between the potentials at its endpoints and the edge parameters.
Coupling edges may therefore transfer additional flow between subinstances, possibly exceeding limits even when each subinstance is feasible on its own.
Here, an additional offset in the potential difference constraints can repair possible infeasibilities of the composite.

\begin{corollary}[Angle compensation in affine potential-based flow problems]\label{cor:potential-join}
Consider network instances $\mathcal{N}_i$ and $\mathcal{N}_j$ on disjoint node sets $V_i$ and $V_j$, with node potentials $\bar\pi^i$ and $\bar\pi^j$ from feasible $\Tract$ solutions (see \Cref{prop:feasibility}).
Let $M$ be a set of coupling edges between the instances.
Each $e=(u,v)\in M$ connects $u\in V_i$ to $v\in V_j$ with $i\ne j$.
Suppose the flow of each edge $e\in M$ is defined by
\[
  f_e=a_e(\pi_u-\pi_v-\phi_e),\qquad a_e>0.
\]
For each such edge, let $\phi_e^{\mathrm{before}}$ be its offset before compensation, taking $\phi_e^{\mathrm{before}}=0$ if none exists, and set
\[
  \psi_e:=\bar\pi^i_u-\bar\pi^j_v-\phi_e^{\mathrm{before}},
  \qquad
  \phi_e:=\phi_e^{\mathrm{before}}+\psi_e.
\]
Suppose that, for every $e\in M$, this offset is admissible and the edge bounds allow $f_e=0$.
If no other constraint introduced with $e$ excludes these values, then the constituent solutions extend to a feasible solution for $\Tract$ on the composed instance after adding any number of such edges.
\end{corollary}
\begin{proof}
Retain the constituent solutions and use the offsets above.
For every $e=(u,v)\in M$ with $u\in V_i$ and $v\in V_j$, we obtain
\begin{align*}
  f_e
  &=a_e\bigl(\bar\pi^i_u-\bar\pi^j_v-(\phi_e^{\mathrm{before}}+\psi_e)\bigr) \\
  &=a_e\bigl(\bar\pi^i_u-\bar\pi^j_v-(\bar\pi^i_u-\bar\pi^j_v)\bigr)=0.
\end{align*}
The zero coupling flows preserve the node balances, and all other existing constraints remain satisfied.
The hypotheses ensure that the offset restrictions, edge bounds, and remaining new constraints are satisfied as well.
\end{proof}

As with variable supply compensation, $\psi_e$ may vary between $0$ (no angle difference compensation) and $\bar\pi^i_u-\bar\pi^j_v-\phi_e^{\mathrm{before}}$ (full compensation, no exchange between networks required), since this range contains the value used in the proof.

Extensions of this corollary to other network flow problems may be possible in similar ways.
Note that the compensation of subinstances to make them self-contained may lead to decreased difficulty of the composite, an issue we examine in \Cref{sec:eval:feasibility}.

\subsection{Security-constrained optimal transmission switching}\label{sec:problem:scots}
As one of two instantiations, we consider SCOTS under the linear DC approximation~\cite{hedman_optimal_2009,donkiewicz2026innoptem}.
SCOTS extends SCOPF by adding topology switching.
The \emph{topological action limit} $\TAL$ bounds the number of transmission branches whose status may be changed.
SCOPF is the special case with $\TAL=0$, in which the network topology is fixed.
The contingency set $\aC$ contains sets $\C$ of failed branches.
The \emph{base case} $\C=\emptyset$ has no branch failures.
A \emph{branch overload variable} is a nonnegative slack variable by which a transmission line may exceed its thermal limit, at a penalty in the objective.
Here, \emph{redispatch cost} is the cost of changing generator outputs from a prescribed reference dispatch.

\problemstatement
  {Scots}
  {A power system graph $G(\mathcal{N})=(V,E)$, fixed loads $\StatD{v}$ and generation bounds $\DemandMin{v},\DemandMax{v}$ for $v\in V$, generation costs, line parameters and operating limits, contingencies $\aC$ including the base case, and a topological action limit $\TAL$.}
  {Generation $\Demand{v}$, voltage angles $\PotC{v}$, line flows $\FlowC{u}{v}$, branch status decisions $\z{u}{v}$ and nonnegative branch overload variables $\FlowOverC{u}{v}{\C}$ satisfying the topological action limit, nodal balance, Kirchhoff voltage laws, generation limits and the relaxed thermal constraints for every $\C\in\aC$.}
  {Minimize redispatch costs and branch overload penalties.}

This specialization has a single commodity, so we omit the commodity index.
The net supply at node $v$ is $s(v)=\Demand{v}-\StatD{v}$.
The generation bounds and fixed demands determine the supply limits in $\mathcal{S}$.
$\mathcal{P}$ contains the thermal limits and the remaining physical and cost parameters.

Here, $\Opt$ is the DC-SCOTS variant stated by \citet{donkiewicz2026innoptem}. 
\todo{Add the arXiv preprint identifier to the InnOpTEM citation when available.}
We refer to that paper for the full mathematical formulation and solution method.
Branch overload is penalized rather than prohibited, so a feasible solution need not be free of overload. 
Our feasibility check $\Tract$ is a DC-OPF (without contingencies) problem with fixed original topology and strict branch limits.
The subinstances are feasible for $\Tract$ by \Cref{prop:feasibility}.
The composite's feasibility for $\Tract$ may be guaranteed by \Cref{cor:potential-join}, since DC power flow problems are potential-based flows with voltage angles $\Pot{v}$ corresponding to the node potentials $\pi_v$.

\subsection{Survivable network design}\label{sec:problem:snd}
Our second instantiation is SND, represented as capacitated modular multicommodity network design~\cite{orlowski2010sndlib,frangioni2009reformulations}, using the formulation of \citet{donkiewicz2026adaptive}.
A \emph{capacity module} is one discrete unit of capacity that can be installed on an edge at a specified cost; several modules may be installed on the same edge.
The contingency set $\aC$ contains subsets $\C\subseteq E$ specifying the failed edges and includes the base case $\emptyset$.

\problemstatement
  {Snd}
  {A graph $G(\mathcal{N})=(V,E)$, commodities $d\in K$ with source $o_d$, sink $t_d$ and demand $b_d>0$, available module types with their capacities and installation costs on each edge $e\in E$, and contingencies $\aC$ including the base case.}
  {Module counts $y_{e,m}\in\Z_{\geq0}$ for each edge $e\in E$ and available module type $m$, and commodity flows $f_{\C,d}$ for $\C\in\aC$ and $d\in K$, satisfying the demands and shared edge capacities, with zero flow on failed edges.}
  {Minimize capacity installation cost.}

For each commodity $d$, the supply limits in $\mathcal{S}$ are fixed to $b_d$ at $o_d$, $-b_d$ at $t_d$, and zero elsewhere.
$\mathcal{P}$ contains the installed capacities and routing costs.
All commodities and both directions share the capacity installed on each edge.

Here, $\Opt$ is SND, while $\Tract$ is a minimum-cost multicommodity flow problem with continuous flows, capacities fixed during preprocessing, and no failures.
After cut flow compensation (\Cref{prop:feasibility}), assigning zero flow to the coupling edges gives a feasible solution for $\Tract$ on the composite.
Coupling edges may carry flow if routing through another source network is cheaper.
To introduce demands between source networks, new demands between nodes of both source networks can be added.
In this case, feasibility must be checked again.

\section{Composition Algorithm}\label{sec:pipeline}

We produce a network instance $\mathcal{N}^{\star}$ in six stages, followed by a feasibility check.
The algorithm draws $s\geq2$ source instances feasible for $\Tract$ from a preprocessed library, extracts a connected subinstance from each, compensates the flow on removed cut edges, selects merge vertices, and adds coupling edges between them to join the subinstances.
We describe the algorithm using the terminology of the abstract network instance of \Cref{sec:problem:abstract}.
Because the merge vertices for each join are selected on the growing composite $C_r$, the joins are performed sequentially.

\paragraph{Source selection.}\label{sec:pipeline:sources}
The library $\mathfrak{N}$ is preprocessed once.
An application may first restrict it using a declared Boolean filter function $h:\mathfrak{N}\to\{0,1\}$, where $h(\mathcal{N})=1$ means that $\mathcal{N}$ is eligible.
For example, $h$ may select supported data features by filtering on $\mathcal{P}$ or a size range by filtering on $|V|$.
By default, we use the Euler bound on the underlying simple graph as a density filter, although violating it does not imply high skewness.
We exclude sources violating the bound before solving $\Tract$ for all remaining instances, and the feasible ones with their solutions form the sampling pool
\[
  \mathfrak{L}
  =\left\{\mathcal{N}\in\mathfrak{N} :
    \begin{aligned}
      &h(\mathcal{N})=1,\ |V|\geq3,\ |E|\leq3|V|-6,\\
      &\Tract(\mathcal{N})\text{ is feasible}
    \end{aligned}
  \right\}.
\]
For each generation attempt, we sample a number of sources $s$ from the integers in $[s_{\min},s_{\max}]$, where $s_{\min}\geq2$.
We then draw the $s$ sources $\mathcal{N}_i$, $i\in\{1,\dots,s\}$, independently and uniformly from $\mathfrak{L}$ with replacement, relabeling if necessary to make all node sets disjoint.
The stored solution provides a flow $f_i(e,k)$ for every edge $e\in E_i$ and commodity $k\in K_i$ for the subsequent cut flow compensation.

\paragraph{Subgraph extraction.}\label{sec:pipeline:subgraph}
For each requested instance $\mathcal{N}^{\star}$, a target size $T$ is drawn uniformly at random from the integers in $[T_{\min},T_{\max}]$.
The parameters $0\leq f_{\min}\leq f_{\max}\leq1$ give lower and upper bounds on the fraction of vertices retained from each source.
Each subinstance size $t_i$ must be an integer in $[\max\{3,\lceil f_{\min}|V_i|\rceil\},\lfloor f_{\max}|V_i|\rfloor]$.
If the selected sources cannot together provide exactly $T$ vertices within these bounds, we redraw the sources.
We then process the sources in random order, drawing each size uniformly from the values that leave a feasible allocation of the remaining target to the remaining sources.
The last source receives the remainder, ensuring $\sum_{i=1}^{s}t_i=T$.
Each instance is then reduced to a connected subgraph with $V_i'\subseteq V_i$ and $|V_i'|=t_i$ by a randomized breadth-first search from a uniformly sampled starting node.
The induced instance $\mathcal{N}_i'$ inherits all node and edge attributes from $\mathcal{N}_i$.
For the later cut flow compensation, the cut edges $\delta(V_i') = \{(u,v)\in E_i : |\{u,v\}\cap V_i'| = 1\}$ are recorded.

\paragraph{Cut flow compensation.}\label{sec:pipeline:boundary}
Using the stored solution $f_i$, we compute the net cut flow $\Delta s_k(v)$ from \eqref{eq:cut-flow-compensation} for each retained node $v$ incident to a cut edge and each commodity $k$.
We adjust the supply bounds using fixed or variable compensation as defined in \Cref{prop:feasibility}, with the choices for each application specified in \Cref{sec:eval:compensation}.
When nodes incident to cut edges also become merge vertices, variable compensation can reduce flows on the coupling edges and nearby edges.
We denote the compensated subinstance by $\mathcal{N}_i'$ and its underlying graph by $G_i=G(\mathcal{N}_i')$.

\paragraph{Merge vertex selection.}\label{sec:pipeline:nodes}
At join $r\in\{2,\dots,s\}$, the current composite $C_{r-1}$ is joined with the next extracted graph $G_r$.
Let
\[
n_r=\min\{|V(C_{r-1})|,|V(G_r)|\}.
\]
We request
\[
q_r=\max\{3,\lceil\log_2 n_r\rceil\}
\]
merge vertex pairs.

We use at least three coupling edges to provide multiple connections between the subinstances.
The logarithmic term gives larger subgraphs more connections with a reasonable growth rate.
Each method produces candidate sets $R_{r-1}\subseteq V(C_{r-1})$ and $R_r\subseteq V(G_r)$.
We select $q_r$ merge vertex pairs $(a_1,b_1),\dots,(a_{q_r},b_{q_r})\in R_{r-1}\times R_r$ so that no vertex occurs in more than one pair, redrawing sources if $|R_{r-1}|<q_r$ or $|R_r|<q_r$.
We compare four merge vertex selection methods:
\begin{itemize}
  \item[a)] \emph{random} (R): sample $q_r$ distinct vertices uniformly from each graph; pair in sampling order;
  \item[b)] \emph{high-degree} (HD): take the $q_r$ highest-degree vertices of each graph; pair in descending order;
  \item[c)] \emph{Kamada--Kawai} (KK): compute a two-dimensional layout whose vertex positions reflect graph distances, take the vertices on the boundary of the convex hull of each layout as candidates, place the layouts side by side and greedily select the shortest cross-graph pairs without reusing a vertex;
  \item[d)] \emph{outer face} (OF): compute and embed a planar subgraph of each graph (heuristically or using branch-and-cut, see \Cref{sec:pipeline:planar}).
  Choose a face with the most distinct incident vertices as the outer face.
  Starting at a random vertex on each outer face, select $q_r$ consecutive vertices in the boundary order from \Cref{sec:problem:composition} and pair them in order.
  By \Cref{lem:boundary-order}, the resulting matching is noncrossing.
\end{itemize}

For a graph with $n=|V|$ vertices and $m=|E|$ edges, random candidates can be selected in $O(n)$ time~\cite{vitter1985random} and high-degree candidates in $O(n+m)$ time~\cite{blum1973time}.
The classical Kamada--Kawai method takes $O(n^3+I_{\mathrm{KK}}n)$ time, where $I_{\mathrm{KK}}$ is the total number of desired updates (improvements) to the vertex positions~\cite[Section~3.3]{kamada1989algorithm}.
Once a planar subgraph has been computed, a planar embedding and the boundary walks of its faces can be computed in $O(n+m)$ time~\cite{boyer2004cutting}.
Computing a maximum planar subgraph, by contrast, is NP-hard~\cite{yannakakis1978nodeedge}.

\paragraph{Sequential composition.}\label{sec:pipeline:merging}
Let $\mathcal{I}_r$ denote the composite instance after joining the first $r$ subinstances, with $\mathcal{I}_1=\mathcal{N}_1'$ and $C_r=G(\mathcal{I}_r)$.
At join $r$, we obtain $\mathcal{I}_r$ by adding the coupling edges $M_r$ between $\mathcal{I}_{r-1}$ and $\mathcal{N}_r'$, one for each selected merge vertex pair.
For the outer face method with maximum planar subgraphs, \Cref{lem:floor,cor:sequential-floor} guarantee that no join adds skewness.
\Cref{lem:suboptimal-mps} bounds the added skewness when the planar subgraphs are not maximum.
The resulting instance $\mathcal{N}^{\star}=\mathcal{I}_s$ inherits the supplies, demands, capacities and other parameters of the extracted instances.
Only the new coupling edges require additional edge parameters, which are specified for each application.

\paragraph{Instance normalization.}\label{sec:pipeline:normalization}
Two source instances may have parameters (e.g., capacity or supply/demand values) measured on different scales.
After the final join, an optional, problem-specific transformation can normalize these attributes without changing the topology or the selected merge vertices.
It should not alter relative parameter values within one network, but only the scale of those values between composed networks.
There is no universal normalization for the abstract model.

\paragraph{Feasibility check.}\label{sec:pipeline:gate}
Although cut flow compensation preserves the feasibility of extracted subinstances under the hypotheses of \Cref{prop:feasibility}, composing feasible subinstances need not preserve feasibility when the hypotheses of neither \Cref{prop:feasibility} nor \Cref{cor:potential-join} hold.
The feasibility check solves $\Tract(\mathcal{N}^{\star})$ after composition, discarding the attempt unless a feasible solution is found.
For each requested output, the algorithm makes at most $A_{\max}$ attempts and reports failure if none succeeds.

\subsection*{Planarization algorithm}\label{sec:pipeline:planar}
We compute planar subgraphs using either a heuristic (OF-H) or branch-and-cut (OF-B) from the Open Graph Drawing Framework (OGDF)~\cite{chimani2014ogdf}, augmented by two additions: valid inequalities derived from the Euler bound, and the reuse of planar subgraphs and skewness bounds in subsequent joins.
In both OF-H and OF-B, planarization is performed separately on the graph's \emph{blocks}, that is, its maximal connected subgraphs that remain connected after deleting any one vertex.
A bridge (an edge whose removal disconnects the graph) forms a block of two vertices.
Planar blocks are retained without deleting edges.
The implementations, solver settings, and resource limits are specified in \Cref{sec:eval:solver-settings}.

\paragraph{Euler inequalities.}
A simple planar graph on $n\geq3$ vertices has at most $3n-6$ edges, or $2n-4$ if it contains no triangle (i.e., $c(n-2)$ with $c=3$ or $c=2$, respectively).
Thus, a subgraph exceeding the applicable bound requires edge deletions.
For $S\subseteq V$ with $|S|\geq3$, let $F$ be the set of edges with both endpoints in $S$.
We use the subgraph Euler inequalities~\cite[Section~3]{chimani2018cycles}
\[
  \sum_{e\in F}x_e\geq |F|-c(|S|-2),
\]
where $x_e\in\{0,1\}$ indicates whether edge $e$ is deleted.
The choice $c=3$ is always valid, and $c=2$ is valid when no triangle has all three vertices in $S$.
To find useful subsets, we adapt the vertex removal procedure of~\cite{charikar2000greedy}: starting with each nonplanar block, we repeatedly remove a vertex of minimum current degree until three remain.
We evaluate the right-hand side on the initial vertex set and after each removal, using $c=3$ and, whenever no triangles remain, also $c=2$.
For each choice of $c$, we retain a subset maximizing the right-hand side $|F|-c(|S|-2)$ along this removal sequence, provided the maximum is positive.
In OF-B, the resulting at most two inequalities per block are added as constraints to the initial formulation, while in OF-H they are used only as lower bounds.

\paragraph{Reusing planar subgraphs and bounds across joins.}\label{sec:pipeline:reuse}
Since we compose the graphs sequentially, we can reuse planar subgraphs and bounds from earlier joins.
Let $H_{r-1}^C$ and $H_r^G$ be the selected planar subgraphs of $C_{r-1}$ and $G_r$.
Their OF join $H_r=(H_{r-1}^C\mathbin{\dot\cup}H_r^G)+M_r$ is planar by \Cref{lem:boundary-order} even if the selected subgraphs are not maximum.
In both OF-B and OF-H, we compare $H_r$ with a fresh heuristic solution and retain the one with more edges.
For OF-H, this is the solution, and for OF-B, it is supplied as an initial feasible solution.

Let $L_{r-1}^C$ and $L_r^G$ be the computed skewness lower bounds for $C_{r-1}$ and $G_r$.
For the next planarization of $C_r$, it still holds that
\[
  \sum_{e\in E(C_{r-1})}x_e\geq L_{r-1}^C,
  \qquad
  \sum_{e\in E(G_r)}x_e\geq L_r^G.
\]
In OF-B, however, since we decompose the problem into blocks, each constraint is added only if its edge set lies entirely within a single block of $C_r$.

Additionally, since $C_{r-1}$ and $G_r$ have disjoint edge sets,
\[
   \operatorname{sk}(C_r) \geq L_{r-1}^C+L_r^G
\]
always holds.

If the available lower bound equals the number of edges omitted by the retained planar subgraph of $C_r$, optimality is proved and no exact solve is needed.

\section{Computational Study}\label{sec:evaluation}\label{sec:eval:sources}
We generate instances from PGLib-OPF v23.07~\cite{pglib} (commit \texttt{dc6be4b}) and SNDlib~\cite{orlowski2010sndlib}.
We use size eligibility filters
\[
\begin{aligned}
  h_{\mathrm{PGLib}}(\mathcal{N})&=
  \begin{cases}
    1, & 73\leq |V(G(\mathcal{N}))|\leq 2000,\\
    0, & \text{otherwise},
  \end{cases}\\[1ex]
  h_{\mathrm{SNDlib}}(\mathcal{N})&=
  \begin{cases}
    1, & 2|K|\,|E(G(\mathcal{N}))|\leq 2{,}000{,}000,\\
    0, & \text{otherwise},
  \end{cases}
\end{aligned}
\]
where $K$ is the commodity set of $\mathcal{N}$ and $2|K|\,|E(G(\mathcal{N}))|$ is the number of flow variables in $\Tract$ for SND.
\Cref{tab:source-screening} summarizes the exclusions and resulting source pools.

\begin{table}[!htbp]
  \scriptsize
  \centering
  \caption{Source screening.}
  \label{tab:source-screening}
  \begin{tabular}{@{}lccccc@{}}
    \toprule
    Library & Initial instances & \multicolumn{3}{c}{Exclusions} & $|\mathfrak{L}|$ \\
    \cmidrule(lr){3-5}
    & & Euler & $\Tract$ & $h$ & \\
    \midrule
    PGLib & 66 & 0 & 1 & 49 & 16 \\
    SNDlib & 26 & 5 & 4 & 1 & 16 \\
    \bottomrule
  \end{tabular}
\end{table}

\paragraph{Source skewness.}

For 28 of the 32 sources, the bounds in \Cref{tab:app-source-skewness,tab:snd-source-skewness} show that deleting at most \SI{5}{\percent} of the edges suffices for planarity.
SNDlib has consistently low skewness, whereas PGLib exhibits a wider range.
Removing branches between different nominal voltage levels makes additional PGLib sources planar.

\begin{table}[!htbp]
  \centering
  \caption{Source skewness intervals (single values are exact) and relative upper bounds in gray. Time limit: \SI{24}{\hour}. PGLib: $s$ with all lines, $s'$ without lines between voltage levels, and $E'$ the remaining edge set.}
  \label{tab:source-skewness}
  \scriptsize
  \setlength{\tabcolsep}{1.8pt}
  \begin{subtable}[t]{0.65\textwidth}
    \centering
    \caption{PGLib}
    \label{tab:app-source-skewness}
    \begin{tabular}{lrrr>{\columncolor{gray!10}}rr>{\columncolor{gray!10}}r}
  \toprule
  Source & $|V|$ & $|E|$ & $[\underline s,\bar s]$ & $\bar s/|E|$ & $[\underline{s}',\bar{s}']$ & $\bar{s}'/|E'|$\\
  \midrule
  \texttt{case73\_ieee\_rts} & 73 & 108 & $4$ & \SI{3.7}{\percent} & $0$ & \SI{0.0}{\percent} \\
  \texttt{case89\_pegase} & 89 & 206 & $[40,\,55]$ & $\leq$\,\SI{26.7}{\percent} & $[31,\,32]$ & $\leq$\,\SI{19.9}{\percent} \\
  \texttt{case118\_ieee} & 118 & 179 & $0$ & \SI{0.0}{\percent} & $0$ & \SI{0.0}{\percent} \\
  \texttt{case162\_ieee\_dtc} & 162 & 280 & $[16,\,30]$ & $\leq$\,\SI{10.7}{\percent} & $0$ & \SI{0.0}{\percent} \\
  \texttt{case179\_goc} & 179 & 222 & $1$ & \SI{0.5}{\percent} & $1$ & \SI{0.5}{\percent} \\
  \texttt{case197\_snem} & 197 & 223 & $0$ & \SI{0.0}{\percent} & $0$ & \SI{0.0}{\percent} \\
  \texttt{case200\_activ} & 200 & 245 & $[1,\,6]$ & $\leq$\,\SI{2.4}{\percent} & $1$ & \SI{0.6}{\percent} \\
  \texttt{case240\_pserc} & 240 & 348 & $10$ & \SI{2.9}{\percent} & $1$ & \SI{0.4}{\percent} \\
  \texttt{case300\_ieee} & 300 & 409 & $9$ & \SI{2.2}{\percent} & $2$ & \SI{0.7}{\percent} \\
  \texttt{case500\_goc} & 500 & 650 & $[14,\,18]$ & $\leq$\,\SI{2.8}{\percent} & $2$ & \SI{0.4}{\percent} \\
  \texttt{case588\_sdet} & 588 & 677 & $[10,\,13]$ & $\leq$\,\SI{1.9}{\percent} & $0$ & \SI{0.0}{\percent} \\
  \texttt{case793\_goc} & 793 & 904 & $[11,\,17]$ & $\leq$\,\SI{1.9}{\percent} & $2$ & \SI{0.3}{\percent} \\
  \texttt{case1354\_pegase} & 1354 & 1710 & $[24,\,50]$ & $\leq$\,\SI{2.9}{\percent} & $[10,\,11]$ & $\leq$\,\SI{0.7}{\percent} \\
  \texttt{case1888\_rte} & 1888 & 2308 & $[25,\,54]$ & $\leq$\,\SI{2.3}{\percent} & $[9,\,10]$ & $\leq$\,\SI{0.6}{\percent} \\
  \texttt{case1951\_rte} & 1951 & 2375 & $[26,\,57]$ & $\leq$\,\SI{2.4}{\percent} & $10$ & \SI{0.5}{\percent} \\
  \texttt{case2000\_goc} & 2000 & 2806 & $[59,\,225]$ & $\leq$\,\SI{8.0}{\percent} & $13$ & \SI{0.6}{\percent} \\
  \bottomrule
\end{tabular}

  \end{subtable}%
  \hfill
  \begin{subtable}[t]{0.33\textwidth}
    \centering
    \caption{SNDlib}
    \label{tab:snd-source-skewness}
    \begin{tabular}{lrrr>{\columncolor{gray!10}}r}
  \toprule
  Source & $|V|$ & $|E|$ & $[\underline s,\bar s]$ & $\bar s/|E|$\\
  \midrule
  \texttt{abilene} & 12 & 15 & $0$ & \SI{0.0}{\percent} \\
  \texttt{atlanta} & 15 & 22 & $0$ & \SI{0.0}{\percent} \\
  \texttt{cost266} & 37 & 57 & $0$ & \SI{0.0}{\percent} \\
  \texttt{france} & 25 & 45 & $0$ & \SI{0.0}{\percent} \\
  \texttt{geant} & 22 & 36 & $2$ & \SI{5.6}{\percent} \\
  \texttt{germany50} & 50 & 88 & $2$ & \SI{2.3}{\percent} \\
  \texttt{giul39} & 39 & 86 & $4$ & \SI{4.7}{\percent} \\
  \texttt{india35} & 35 & 80 & $4$ & \SI{5.0}{\percent} \\
  \texttt{janos-us-ca} & 39 & 61 & $0$ & \SI{0.0}{\percent} \\
  \texttt{janos-us} & 26 & 42 & $0$ & \SI{0.0}{\percent} \\
  \texttt{nobel-eu} & 28 & 41 & $0$ & \SI{0.0}{\percent} \\
  \texttt{nobel-germany} & 17 & 26 & $0$ & \SI{0.0}{\percent} \\
  \texttt{nobel-us} & 14 & 21 & $1$ & \SI{4.8}{\percent} \\
  \texttt{norway} & 27 & 51 & $1$ & \SI{2.0}{\percent} \\
  \texttt{sun} & 27 & 51 & $1$ & \SI{2.0}{\percent} \\
  \texttt{zib54} & 54 & 80 & $0$ & \SI{0.0}{\percent} \\
  \bottomrule
\end{tabular}

  \end{subtable}
\end{table}

\subsection{Experimental setup}\label{sec:eval:setup}\label{sec:eval:compensation}

\paragraph{Generation settings.}
For both applications, we compare five methods on the same 100 accepted inputs: random (R), high-degree (HD), Kamada--Kawai (KK), and outer face selection with either heuristic (OF-H) or branch-and-cut (OF-B) planarization.
For each input, the extracted graphs, their order, and all random seeds are fixed across methods.
We accept an input only if all five composites are feasible for $\Tract$.
\Cref{tab:generation-settings} lists the generation settings.

\begin{table}[!htbp]
  \centering
  \scriptsize
  \caption{Generation settings.}
  \label{tab:generation-settings}
  \begin{tabular}{@{}lccc@{}}
    \toprule
    Setting & Symbols & PGLib & SNDlib \\
    \midrule
    Source networks & $[s_{\min},s_{\max}]$ & $[2,5]$ & $[2,5]$ \\
    Retained vertex fraction & $[f_{\min},f_{\max}]$ & $[0.15,0.85]$ & $[0.15,0.85]$ \\
    Target vertices & $[T_{\min},T_{\max}]$ & $[100,1000]$ & $[12,100]$ \\
    Maximum attempts & $A_{\max}$ & 25 & 25 \\
    Outer face vertex selection & & Consecutive & Consecutive \\
    Cut flow compensation & & Variable & Fixed \\
    \bottomrule
  \end{tabular}
\end{table}

For PGLib, the variable cut flow compensation is implemented using auxiliary generators with zero cost and additional loads.
We refrain from using the additional voltage angle offsets proposed in \Cref{cor:potential-join}.
We define coupling edges as transformers with reactance $10^{-4}$ per unit, unlimited thermal capacity, tap ratio 1, and zero phase shift.

For SNDlib, we decompose flows from $\Tract$ into paths terminating at the original terminals or, if one of those was removed during subgraph extraction, at the boundary vertex they flowed through.
Coupling edges have zero module, routing, and setup costs.
Their module capacity equals the combined positive demand of the two subinstances.

\paragraph{Computational settings and measurements.}
\label{sec:eval:solver-settings}
We use OGDF version foxglove-202510, with the additions described in \Cref{sec:pipeline:planar}.
OF-H uses OGDF's heuristic mode (\texttt{PlanarSubgraphFast} with its default ten runs).
OF-B uses its exact mode (\texttt{MaxCPlanarMaster} with default settings), with the ABACUS~\cite{junger1998abacus} branch-and-cut solver and COIN-OR's linear programming solver~\cite{clpguide}.

For PGLib, both $\Tract$ and $\Opt$ are solved using the exact SCOTS solver of \citet{donkiewicz2026innoptem} with Gurobi~13.0.2.
For SNDlib, we solve $\Tract$ using our own implementation with Gurobi~13.0.3 and $\Opt$ using the SND implementation of \citet{donkiewicz2026adaptive} with Gurobi~12.0.3.

We use compute nodes with eight physical cores (Intel Xeon L5630, \SI{2.13}{\giga\hertz}) and at least \SI{16}{GB} of memory.
We use up to sixteen threads for heuristic planarization and one for branch-and-cut, Euler preprocessing, and application solves.
Exact planarization is limited to \SI{600}{\second} per graph during generation and \SI{3600}{\second} for final evaluation, including initialization.
Reported runtimes additionally include Euler preprocessing, solution and bound reuse, and the separate heuristic call before the exact solve in OF-B, all of which do not count towards the above limits.

For the final skewness evaluation, all five methods use the same known planar subgraphs of the original extracted subinstances and the same random seed, and also receive all applicable bounds from generation.
To obtain an initial feasible solution, we join these subgraphs with all coupling edges if the result is planar, and otherwise with a minimal subset that connects the subinstances.

For nonnegative values $z_1,\dots,z_N$, we report the shifted geometric mean
\[
  \exp\!\left(\frac{1}{N}\sum_{i=1}^{N}\log(1+z_i)\right)-1.
\]
We use this mean for runtimes and skewness, and medians for structural summaries.

\FloatBarrier
\subsection{Effects of composition}\label{sec:eval:composite}

\paragraph{Added skewness.}

\begin{figure}[!htbp]
  \centering
  \input{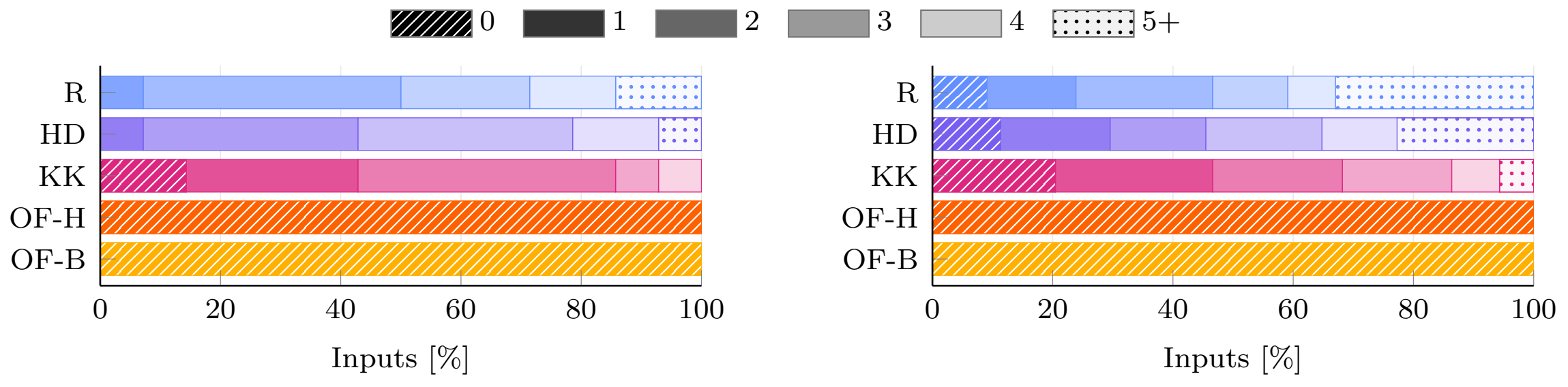}
  \caption{Distribution of added skewness on inputs with exact added skewness values for all five methods. Within each method color, lighter shades indicate larger added skewness.}
  \label{fig:eval-exact-skewness-distribution}
\end{figure}

We first consider the 14 PGLib and 88 SNDlib inputs for which added skewness is known exactly for all five methods.
Both outer face variants add no skewness on these inputs (\Cref{fig:eval-exact-skewness-distribution}).
KK has more cases with zero added skewness than R or HD, but most of its composites still add skewness.
Thus, selecting vertices on the boundary of a geometric layout does not by itself prevent added skewness.

\FloatBarrier

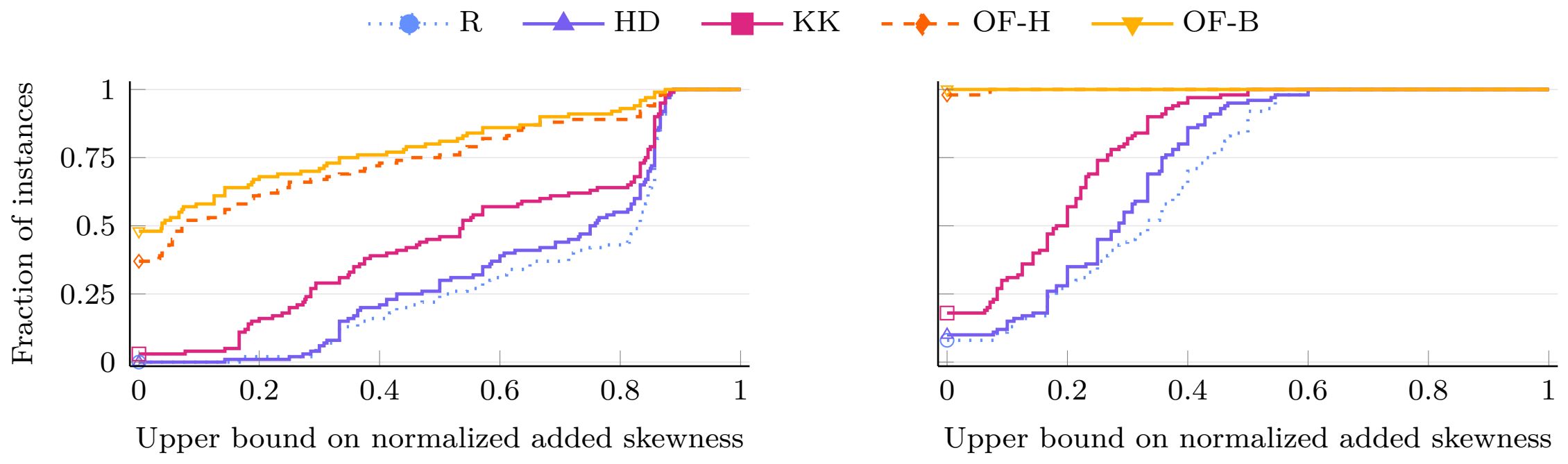
\begin{figure}[!htbp]
  \centering
  \begingroup
    \pgfplotslegendfromname{skewnesslegend}\par\smallskip
    \begin{subfigure}[t]{.485\linewidth}
      \centering
      \begin{tikzpicture}
\begin{axis}[
width=.97\linewidth,height=.5567\linewidth,
axis lines*=left,axis line style={line width=.5pt},
xmin=-.015,xmax=1.015,ymin=-.015,ymax=1.03,
xtick={0,.2,.4,.6,.8,1},ytick={0,.25,.5,.75,1},
scaled ticks=false,tick label style={font=\scriptsize},
label style={font=\scriptsize},
xlabel={Upper bound on normalized added skewness},
ymajorgrids=true,grid style={gray!20,line width=.3pt},
legend columns=5,legend cell align=left,
legend style={font=\scriptsize,draw=none,
/tikz/every even column/.append style={column sep=3mm}},
ylabel={Fraction of instances},
legend to name=skewnesslegend,
]
\addplot+[no marks,line width=.9pt,color=pigMethod0,dash pattern=on .6pt off 2.4pt,forget plot] table[x=x,y=y] {resources/evaluation/plot-data/skewness-pglib-random.dat};
\addplot[only marks,forget plot,mark=o,mark size=1.8pt,color=pigMethod0] coordinates {(0,0)};
\addplot+[no marks,line width=.9pt,color=pigMethod1,solid,forget plot] table[x=x,y=y] {resources/evaluation/plot-data/skewness-pglib-high_degree.dat};
\addplot[only marks,forget plot,mark=triangle,mark size=1.8pt,color=pigMethod1] coordinates {(0,0)};
\addplot+[no marks,line width=.9pt,color=pigMethod2,solid,forget plot] table[x=x,y=y] {resources/evaluation/plot-data/skewness-pglib-kamada_kawai.dat};
\addplot[only marks,forget plot,mark=square,mark size=1.8pt,color=pigMethod2] coordinates {(0,0.03)};
\addplot+[no marks,line width=.9pt,color=pigMethod3,dashed,forget plot] table[x=x,y=y] {resources/evaluation/plot-data/skewness-pglib-heuristic_outer_face.dat};
\addplot[only marks,forget plot,mark=diamond,mark size=1.8pt,color=pigMethod3] coordinates {(0,0.37)};
\addplot+[no marks,line width=.9pt,color=pigMethod4,solid,forget plot] table[x=x,y=y] {resources/evaluation/plot-data/skewness-pglib-exact_outer_face.dat};
\addplot[only marks,forget plot,mark=triangle,mark size=1.8pt,color=pigMethod4,mark options={rotate=180}] coordinates {(0,0.48)};
\addlegendimage{line width=.9pt,color=pigMethod0,dash pattern=on .6pt off 2.4pt,mark=*}
\addlegendentry{R}
\addlegendimage{line width=.9pt,color=pigMethod1,solid,mark=triangle*}
\addlegendentry{HD}
\addlegendimage{line width=.9pt,color=pigMethod2,solid,mark=square*}
\addlegendentry{KK}
\addlegendimage{line width=.9pt,color=pigMethod3,dashed,mark=diamond*}
\addlegendentry{OF-H}
\addlegendimage{line width=.9pt,color=pigMethod4,solid,mark=triangle*,mark options={rotate=180}}
\addlegendentry{OF-B}
\end{axis}
\end{tikzpicture}
      \caption{PGLib}\label{fig:eval-skewness-pglib}
    \end{subfigure}\hfill
    \begin{subfigure}[t]{.485\linewidth}
      \centering
      \begin{tikzpicture}
\begin{axis}[
width=.97\linewidth,height=.5567\linewidth,
axis lines*=left,axis line style={line width=.5pt},
xmin=-.015,xmax=1.015,ymin=-.015,ymax=1.03,
xtick={0,.2,.4,.6,.8,1},ytick={0,.25,.5,.75,1},
scaled ticks=false,tick label style={font=\scriptsize},
label style={font=\scriptsize},
xlabel={Upper bound on normalized added skewness},
ymajorgrids=true,grid style={gray!20,line width=.3pt},
legend columns=5,legend cell align=left,
legend style={font=\scriptsize,draw=none,
/tikz/every even column/.append style={column sep=3mm}},
yticklabels={},
]
\addplot+[no marks,line width=.9pt,color=pigMethod0,dash pattern=on .6pt off 2.4pt,forget plot] table[x=x,y=y] {resources/evaluation/plot-data/skewness-snd-random.dat};
\addplot[only marks,forget plot,mark=o,mark size=1.8pt,color=pigMethod0] coordinates {(0,0.08)};
\addplot+[no marks,line width=.9pt,color=pigMethod1,solid,forget plot] table[x=x,y=y] {resources/evaluation/plot-data/skewness-snd-high_degree.dat};
\addplot[only marks,forget plot,mark=triangle,mark size=1.8pt,color=pigMethod1] coordinates {(0,0.1)};
\addplot+[no marks,line width=.9pt,color=pigMethod2,solid,forget plot] table[x=x,y=y] {resources/evaluation/plot-data/skewness-snd-kamada_kawai.dat};
\addplot[only marks,forget plot,mark=square,mark size=1.8pt,color=pigMethod2] coordinates {(0,0.18)};
\addplot+[no marks,line width=.9pt,color=pigMethod3,dashed,forget plot] table[x=x,y=y] {resources/evaluation/plot-data/skewness-snd-heuristic_outer_face.dat};
\addplot[only marks,forget plot,mark=diamond,mark size=1.8pt,color=pigMethod3] coordinates {(0,0.98)};
\addplot+[no marks,line width=.9pt,color=pigMethod4,solid,forget plot] table[x=x,y=y] {resources/evaluation/plot-data/skewness-snd-exact_outer_face.dat};
\addplot[only marks,forget plot,mark=triangle,mark size=1.8pt,color=pigMethod4,mark options={rotate=180}] coordinates {(0,1)};
\end{axis}
\end{tikzpicture}
      \caption{SNDlib}\label{fig:eval-skewness-snd}
    \end{subfigure}
  \endgroup
  \caption{Distribution of upper bounds on normalized added skewness on all 100 final instances.}
  \label{fig:eval-skewness}
\end{figure}

To include all 100 inputs per application, we compare upper bounds on added skewness relative to the number of coupling edges (\Cref{fig:eval-skewness} and \Cref{tab:eval-skewness}).
Both outer face methods give substantially smaller upper bounds for PGLib.
With OF-H and OF-B, the upper bound on added skewness is zero for 37 and 48 PGLib inputs and 98 and 100 SNDlib inputs, respectively, compared to at most three and 18 for any of the alternatives.
The bounds show that OF-B adds strictly less skewness than R, HD, and KK on 45 PGLib inputs and 77 SNDlib inputs.
Relative to OF-H, OF-B tightens the upper bound on 37 PGLib inputs and certifies strictly smaller added skewness on eight PGLib and two SNDlib inputs.
\Cref{fig:app-example-composites} shows five composites from the same four subinstances: OF-B preserves their total skewness of three, OF-H adds one, and R, HD, and KK each add at least three.

\FloatBarrier

\newcommand{\CompositeSourceLegend}{\normalsize\begin{tabular}{ll}
\textcolor[HTML]{5B7FA3}{$\bullet$} $G_1$ & \texttt{case179\_goc}\\
\textcolor[HTML]{C98B5A}{$\bullet$} $G_2$ & \texttt{case300\_ieee}\\
\textcolor[HTML]{5E9B76}{$\bullet$} $G_3$ & \texttt{case793\_goc}\\
\textcolor[HTML]{A77C9E}{$\bullet$} $G_4$ & \texttt{case118\_ieee}\\
\end{tabular}}
\newcommand{\CompositeGraphPanel}[1]{\includegraphics[width=\linewidth,height=6.95cm,keepaspectratio]{#1}}
\newcommand{\CompositePanel}[1]{\CompositeGraphPanel{#1}\par\CompositeSourceLegend}

\begin{figure}[p]
  \captionsetup[subfigure]{skip=2pt}
  \centering
  \begin{subfigure}[t]{0.49\textwidth}
    \centering
    \CompositeGraphPanel{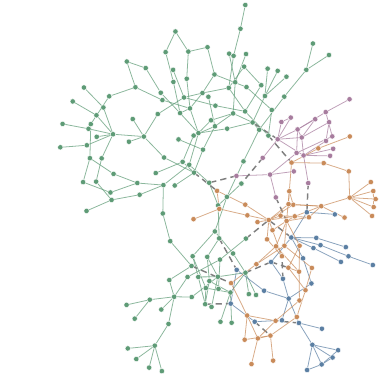}
    \caption{Random: $\operatorname{sk}(C)\in[8,10]$.}
  \end{subfigure}
  \hfill
  \begin{subfigure}[t]{0.49\textwidth}
    \centering
    \CompositeGraphPanel{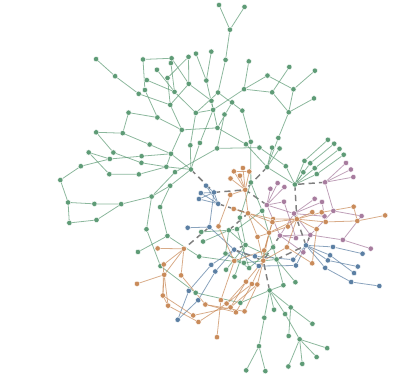}
    \caption{High degree: $\operatorname{sk}(C)\in[8,10]$.}
  \end{subfigure}
  \par\vspace{2mm}
  \begin{subfigure}[c]{0.49\textwidth}
    \centering
    \CompositeGraphPanel{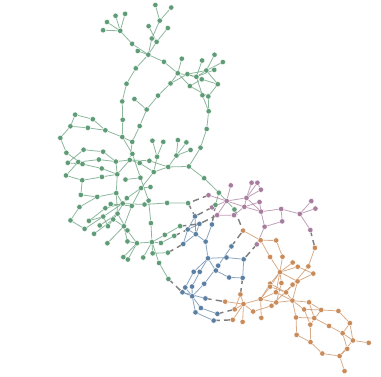}
    \caption{Kamada--Kawai: $\operatorname{sk}(C)=6$.}
  \end{subfigure}
  \hfill
  \begin{minipage}[c]{0.49\textwidth}
    \centering
    \CompositeSourceLegend
  \end{minipage}
  \par\vspace{2mm}
  \begin{subfigure}[t]{0.49\textwidth}
    \centering
    \CompositeGraphPanel{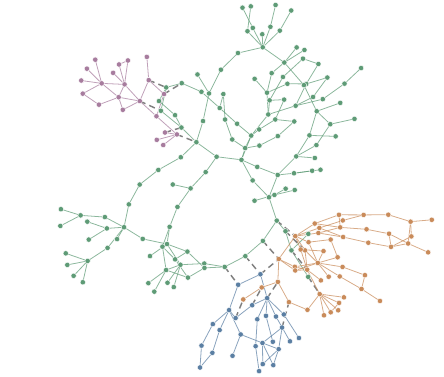}
    \caption{Outer face (heuristic): $\operatorname{sk}(C)=4$.}
  \end{subfigure}
  \hfill
  \begin{subfigure}[t]{0.49\textwidth}
    \centering
    \CompositeGraphPanel{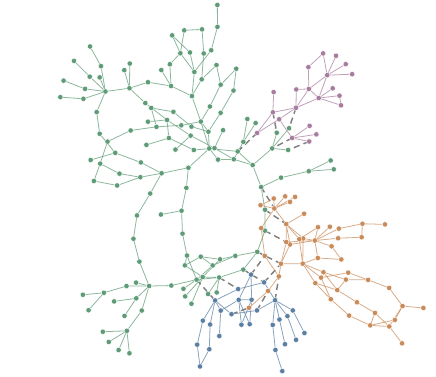}
    \caption{Outer face (branch-and-cut): $\operatorname{sk}(C)=3$.}
  \end{subfigure}
  \caption{Final composites for input \texttt{i070}, drawn with Kamada--Kawai layouts that may visually favor KK. Apparent crossings do not measure skewness.}
  \label{fig:app-example-composites}
\end{figure}

\paragraph{Other structural characteristics.}

\Cref{fig:eval-bounded} compares changes in the degree distribution, mean clustering, and maximum degree with the bounds from \Cref{prop:preserved}.
For the degree distribution, we report the largest absolute change among the fractions of vertices in each degree class.
HD attains the upper bound on the increase in maximum degree on about two thirds of the inputs, which is consistent with repeatedly selecting vertices of maximum degree.
KK tends to select vertices of low degree (\Cref{fig:app-merge-vertex-degrees}), leaving the maximum degree almost unchanged.
Its concentration on vertices of the same degree, particularly degree one in PGLib, explains the largest median change in degree frequencies relative to the bound.
OF and R choose vertices across the degree distribution, yielding smaller median changes in degree frequencies than KK and smaller median increases in maximum degree than HD.

For all methods, median changes in mean clustering are small relative to the bound, indicating that the constituents' local clustering changes little on average, even among the merge vertices.

\Cref{fig:eval-structure} compares diameter, algebraic connectivity, and degree assortativity.
HD and R give smaller median diameters than KK and OF on the same inputs, suggesting that their choice of merge vertices provides shorter paths between subinstances.

For PGLib, OF gives median algebraic connectivity closer to the source median than R and HD, whose medians lie above it.
For SNDlib, all methods give medians below the source median, with the lowest values for OF suggesting slightly more pronounced bottlenecks between subinstances.

HD gives higher median degree assortativity than the other methods, particularly for PGLib, consistent with adding edges between high-degree vertices.
OF keeps median degree assortativity within the source interquartile range, avoiding the pronounced shift in degree mixing observed for HD in PGLib and giving distributions similar to R and KK.

\FloatBarrier
\subsection{Generation and skewness evaluation time}

\begin{figure}[!htbp]
  \centering
  \input{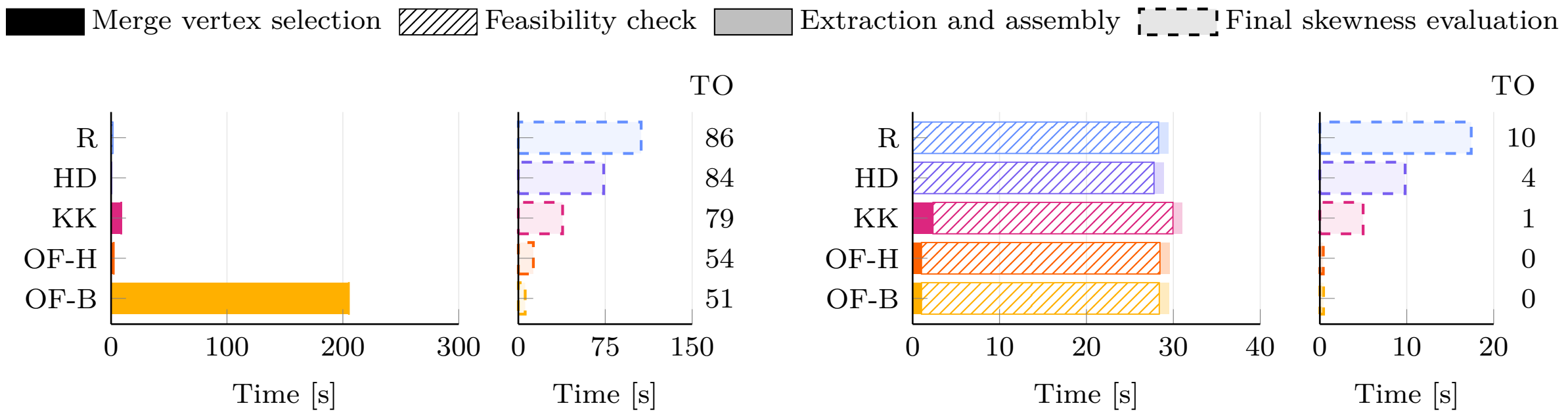}
  \caption{Shifted geometric means; 100 accepted inputs per method. Final skewness times include only exact cases. \textsc{TO}: runs reaching \SI{3600}{\second} in the skewness calculation on the composite.}
  \label{fig:eval-runtime}
\end{figure}

Merge vertex selection dominates the generation time for OF on PGLib, where the higher skewness makes planarization more expensive (\Cref{fig:eval-runtime}).
SNDlib instances have very low skewness and are cheaper to planarize, so feasibility checking for $\Tract$ accounts for most of the generation time.

\Needspace{4\baselineskip}
OF-B takes roughly one hundred times the selection time of OF-H on PGLib, but little additional time on SNDlib.
These results suggest that OF-H is often sufficient, while OF-B offers further improvements at little additional cost on instances as close to planar as those in SNDlib.

Final skewness evaluation reaches the time limit on most PGLib composites, even for OF, but certifies exact skewness for most SNDlib composites.

\subsection{Optimization experiments}\label{sec:eval:feasibility}

We solve $\Opt$ on the generated instances for both applications, using the problem definitions in \Cref{sec:problem:scots,sec:problem:snd}.

For SCOTS, we set $\TAL=1$, permit overload slack, and exclude single-line contingencies that disconnect the network.
SCOPF warmstarting is enabled.
For SND, the contingency set $\aC$ contains the base case and all single-edge failures under which all demand terminals remain in one connected component.
Nodes without demands may become disconnected.
Both applications use a \SI{10800}{\second} solver limit.

\paragraph{Feasibility.}\label{sec:eval:sensitivity}

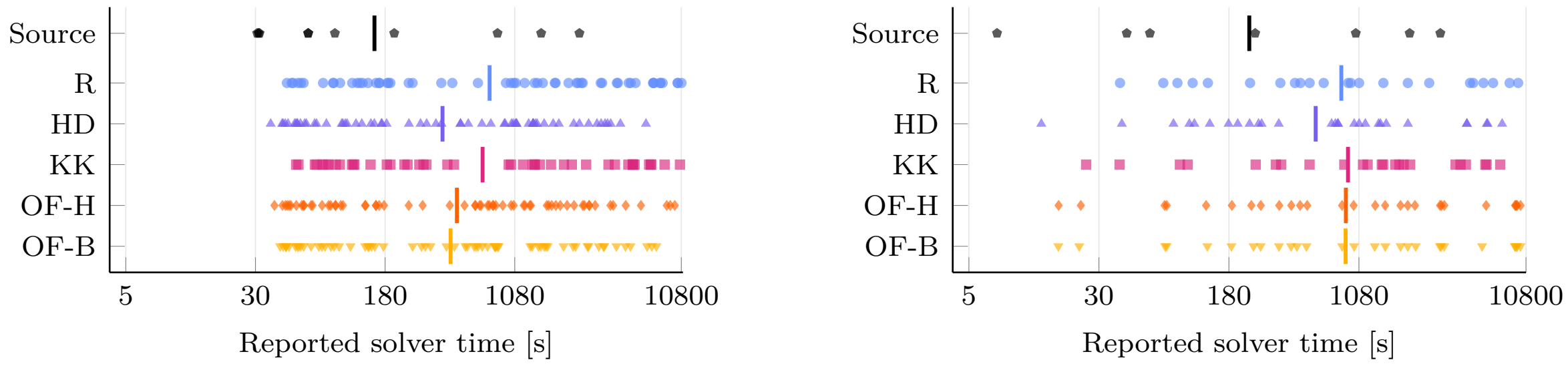
\begin{figure}[!htbp]
  \centering
  \begin{subfigure}[t]{0.475\linewidth}
    \centering
\begingroup
\pigSetRowSpacing
\begin{tikzpicture}
\begin{axis}[
  width=.96\linewidth,height=\pigSourcePlotHeight,
  xmode=log,xmin=4,xmax=11043.1041489,
  xtick={5,30,180,1080,10800},
  xticklabels={5,30,180,1080,10800},
  ymin=\pigSourceYMin,ymax=4.64,ytick={\pigSourceY,0,1,2,3,4},y dir=reverse,
  yticklabels={Source,R,HD,KK,OF-H,OF-B},
  axis lines*=left,axis line style={line width=.5pt},
  tick label style={font=\scriptsize},label style={font=\scriptsize},
  xlabel={Reported solver time [s]},
  xmajorgrids=true,grid style={gray!20,line width=.3pt},
]
\addplot[only marks,mark=*,mark size=1.25pt,draw=pigMethod0,fill=pigMethod0,opacity=.65] table[x=runtime,y=y] {resources/evaluation/plot-data/application-runtime-scopf-random.dat};
\addplot[only marks,mark=|,mark size=5pt,line width=1.1pt,color=pigMethod0] coordinates {(763.22981756,0)};
\addplot[only marks,mark=triangle*,mark size=1.25pt,draw=pigMethod1,fill=pigMethod1,opacity=.65] table[x=runtime,y=y] {resources/evaluation/plot-data/application-runtime-scopf-high_degree.dat};
\addplot[only marks,mark=|,mark size=5pt,line width=1.1pt,color=pigMethod1] coordinates {(398.085943403,1)};
\addplot[only marks,mark=square*,mark size=1.25pt,draw=pigMethod2,fill=pigMethod2,opacity=.65] table[x=runtime,y=y] {resources/evaluation/plot-data/application-runtime-scopf-kamada_kawai.dat};
\addplot[only marks,mark=|,mark size=5pt,line width=1.1pt,color=pigMethod2] coordinates {(692.659995287,2)};
\addplot[only marks,mark=diamond*,mark size=1.25pt,draw=pigMethod3,fill=pigMethod3,opacity=.65] table[x=runtime,y=y] {resources/evaluation/plot-data/application-runtime-scopf-heuristic_outer_face.dat};
\addplot[only marks,mark=|,mark size=5pt,line width=1.1pt,color=pigMethod3] coordinates {(485.480141529,3)};
\addplot[only marks,mark=triangle*,mark options={rotate=180},mark size=1.25pt,draw=pigMethod4,fill=pigMethod4,opacity=.65] table[x=runtime,y=y] {resources/evaluation/plot-data/application-runtime-scopf-exact_outer_face.dat};
\addplot[only marks,mark=|,mark size=5pt,line width=1.1pt,color=pigMethod4] coordinates {(445.176600544,4)};
\addplot[only marks,mark=pentagon*,mark size=1.25pt,draw=pigMethod5,fill=pigMethod5,opacity=.65] table[x=runtime,y expr=\pigSourceY] {resources/evaluation/plot-data/application-runtime-scopf-source.dat};
\addplot[only marks,mark=|,mark size=5pt,line width=1.1pt,color=pigMethod5] coordinates {(155.429768838,\pigSourceY)};
\end{axis}
\end{tikzpicture}
\endgroup
    \caption{SCOTS}
  \end{subfigure}
  \hfill
  \begin{subfigure}[t]{0.475\linewidth}
    \centering
\begingroup
\pigSetRowSpacing
\begin{tikzpicture}
\begin{axis}[
  width=.96\linewidth,height=\pigSourcePlotHeight,
  xmode=log,xmin=4,xmax=10800,
  xtick={5,30,180,1080,10800},
  xticklabels={5,30,180,1080,10800},
  ymin=\pigSourceYMin,ymax=4.64,ytick={\pigSourceY,0,1,2,3,4},y dir=reverse,
  yticklabels={Source,R,HD,KK,OF-H,OF-B},
  axis lines*=left,axis line style={line width=.5pt},
  tick label style={font=\scriptsize},label style={font=\scriptsize},
  xlabel={Reported solver time [s]},
  xmajorgrids=true,grid style={gray!20,line width=.3pt},
]
\addplot[only marks,mark=*,mark size=1.25pt,draw=pigMethod0,fill=pigMethod0,opacity=.65] table[x=runtime,y=y] {resources/evaluation/plot-data/application-runtime-snd-random.dat};
\addplot[only marks,mark=|,mark size=5pt,line width=1.1pt,color=pigMethod0] coordinates {(847.3003255,0)};
\addplot[only marks,mark=triangle*,mark size=1.25pt,draw=pigMethod1,fill=pigMethod1,opacity=.65] table[x=runtime,y=y] {resources/evaluation/plot-data/application-runtime-snd-high_degree.dat};
\addplot[only marks,mark=|,mark size=5pt,line width=1.1pt,color=pigMethod1] coordinates {(594.962571727,1)};
\addplot[only marks,mark=square*,mark size=1.25pt,draw=pigMethod2,fill=pigMethod2,opacity=.65] table[x=runtime,y=y] {resources/evaluation/plot-data/application-runtime-snd-kamada_kawai.dat};
\addplot[only marks,mark=|,mark size=5pt,line width=1.1pt,color=pigMethod2] coordinates {(928.783152944,2)};
\addplot[only marks,mark=diamond*,mark size=1.25pt,draw=pigMethod3,fill=pigMethod3,opacity=.65] table[x=runtime,y=y] {resources/evaluation/plot-data/application-runtime-snd-heuristic_outer_face.dat};
\addplot[only marks,mark=|,mark size=5pt,line width=1.1pt,color=pigMethod3] coordinates {(902.173949917,3)};
\addplot[only marks,mark=triangle*,mark options={rotate=180},mark size=1.25pt,draw=pigMethod4,fill=pigMethod4,opacity=.65] table[x=runtime,y=y] {resources/evaluation/plot-data/application-runtime-snd-exact_outer_face.dat};
\addplot[only marks,mark=|,mark size=5pt,line width=1.1pt,color=pigMethod4] coordinates {(899.179264714,4)};
\addplot[only marks,mark=pentagon*,mark size=1.25pt,draw=pigMethod5,fill=pigMethod5,opacity=.65] table[x=runtime,y expr=\pigSourceY] {resources/evaluation/plot-data/application-runtime-snd-source.dat};
\addplot[only marks,mark=|,mark size=5pt,line width=1.1pt,color=pigMethod5] coordinates {(238.318751765,\pigSourceY)};
\end{axis}
\end{tikzpicture}
\endgroup
    \caption{SND}
  \end{subfigure}
  \caption{Globally optimal runs only. SCOTS times include SCOPF warmstarting. Points: runs; ticks: shifted geometric means.}
  \label{fig:eval-application-runtime}
\end{figure}

We find feasible solutions for \SI{65.6}{\percent} of SCOTS and \SI{42.8}{\percent} of SND composites.
In both domains, we find feasible solutions most often with HD and least often with KK.
The ``source'' comparisons in \Cref{fig:eval-application-outcomes,fig:eval-application-runtime} use the sets in \Cref{tab:source-skewness}.

\Needspace{4\baselineskip}
OF selects merge vertices consecutively along the outer face boundary.
For OF-B on PGLib, we test three alternatives: selecting vertices at equally spaced positions along the boundary, skipping one vertex between selections where possible, and preferring vertices of degree greater than one during consecutive selection.
We test whether spreading out the connections and avoiding vertices of degree one help us find feasible SCOPF solutions.
The tested changes offer little improvement in the number of runs returning a feasible solution (\Cref{app:eval-spacing}).

We also test normalization across voltage levels for all five methods on PGLib.
It does not increase the number of SCOPF runs returning a feasible solution (\Cref{app:eval-normalization}).

\FloatBarrier
\paragraph{Interaction between subinstances.}

We introduce no demands between SND subinstances and compensate flows over cut edges in PGLib, which raises the question of whether the subinstances interact after composition.
We therefore examine whether coupling affects power flows and capacity installation costs.

For PGLib, we examine flows on coupling edges in SCOPF solutions (i.e., SCOTS without any switching), observing that every coupling edge carries nonzero base-case flow.
We also compute the ratio of the mean absolute flow on coupling edges to that on inherited edges within each composite (\Cref{fig:eval-flows-pglib}).
The median ratio is \SI{99.4}{\percent} for OF-H and \SI{104.9}{\percent} for OF-B, with \SI{134.5}{\percent} for R, meaning that the coupling edges transfer roughly as much power on average as other edges in the constituents do.
KK and HD have median ratios of \SI{70.2}{\percent} and \SI{280.1}{\percent}, respectively, indicating lower mean absolute flows on coupling edges than on inherited edges for KK and substantially higher flows for HD.

For SNDlib, routing with fixed capacities and no failures yields nonzero flow between subinstances on all OF instances (\Cref{app:eval-flows}).
After summing flows over all commodities and both directions, the median ratio of the mean flow on added coupling edges to that on inherited edges is \SI{11.1}{\percent} for OF-H and \SI{10.5}{\percent} for OF-B.

For the full SND problem, we also compare the composite objective with the sum of the optimal objectives of the compensated subinstances.
Without flow between subinstances, we would expect these values to coincide under the same contingency requirements.
Among the 130 comparisons where all subinstances and the composite are solved to optimality, the composite objective is lower in 83 cases, equal in one, and higher in 46.
The lower costs indicate that coupling edges carry flow between subinstances, allowing demands to be routed through other subinstances even though no demands connect them.
Higher costs can arise when coupling keeps demand terminals connected after previously excluded edge failures, adding these failures to the contingency set.

We conclude that the subinstances interact after composition in both applications.
For PGLib, OF yields comparable mean absolute flows on coupling and inherited edges.
For SNDlib, coupling can reduce capacity installation costs even without demands between subinstances.
Onward routing could induce additional flow across the join by placing the demand terminals introduced by cut flow compensation farther inside another subinstance.

\FloatBarrier

\paragraph{Computational difficulty.}

Among the generated instances solved to optimality, runtimes range from seconds to about three hours in both applications (\Cref{fig:eval-application-runtime}).
Composition therefore provides benchmark instances at different levels of computational difficulty from the same small source pool.
These results suggest that merge vertex selection affects both feasibility and computational difficulty.

\FloatBarrier

\section{Summary}\label{sec:conclusion}

Meaningful computational studies require sufficiently large and diverse test sets, but the availability of suitable real-world network instances is often limited.
We generate additional infrastructure network instances by repeatedly extracting and joining connected subinstances in different ways, reusing their topology and parameters.
We derive bounds on added skewness and structural changes, and propose a strategy for selecting the vertices at which subinstances are joined.
This strategy adds no skewness when maximum planar subgraphs are computed.
To reduce the number of infeasible instances, we compensate flows over cut edges and solve a tractable auxiliary problem during generation.

In experiments on PGLib-OPF and SNDlib, two libraries of infrastructure networks, our outer face selection method gives the smallest mean upper bounds on added skewness.
Using branch-and-cut for planarization, our method certifies no added skewness for 48 of 100 PGLib inputs and all 100 SNDlib inputs.
The optimization results show that the joined subinstances interact even after cut flow compensation, with coupling affecting power flows and capacity installation costs.
Among the generated instances solved to optimality, SCOTS and SND runtimes range from seconds to about three hours, demonstrating variation in computational difficulty.

Our method can thus generate many additional test instances from a few existing networks, or even a single one.
Reusing their topology and parameters reduces the need for domain knowledge compared with generating realistic networks from scratch.
Controlling added skewness also helps preserve near-planarity, a characteristic of realistic infrastructure networks.

\paragraph*{Author Contributions.}\label{sec:contributions}
\paperAuthor: Conceptualization, Data curation, Formal analysis, Investigation, Methodology, Software, Validation, Visualization, Writing -- original draft, Writing -- review \& editing.

\paragraph*{Funding.}\label{sec:funding}
This work received no specific funding.

\paragraph*{Competing Interests.}\label{sec:interests}
The author declares no competing interests.

\paragraph*{Use of Large Language Models.}\label{sec:llm}
OpenAI Codex and Anthropic Claude assisted with code development, theoretical discussions and manuscript revisions.
The author reviewed and edited all outputs and takes full responsibility for the work.

\paragraph*{Data and Code Availability.}\label{sec:availability}
The source code, generated test instances, solver settings, and reproduction instructions are available in \href{https://github.com/tidonk/PlanarInstanceGeneration/tree/paper}{PlanarInstanceGeneration}, tag \texttt{paper}.

\bibliography{literature}

\paperappendix
\clearpage

\appendix

\section{Additional Evaluation Results}\label{app:figures}
\subsection{Skewness bounds}\label{app:eval-skewness}

For $p$ extracted subgraphs, we measure added skewness and normalize it by the number of coupling edges:
\[
  \begin{aligned}
    a(C)&=\operatorname{sk}(C)-\sum_{i=1}^{p}\operatorname{sk}(G_i),\\
    \eta(C)&=\frac{a(C)}{|M|}.
  \end{aligned}
\]
Let $L_Q$ and $U_Q$ denote the computed lower and upper bounds on $\operatorname{sk}(Q)$, respectively.
By \Cref{lem:skewness-bounds}, bounds on added skewness are
\[
  \begin{aligned}
    \underline a&=\max\!\left\{0,L_C-\sum_{i=1}^{p}U_{G_i}\right\},\\
    \bar a&=\min\!\left\{|M|-(p-1),U_C-\sum_{i=1}^{p}L_{G_i}\right\}.
  \end{aligned}
\]
Dividing these bounds by $|M|$ gives $\underline\eta$ and $\bar\eta$.
Thus, $[\underline a,\bar a]$ bounds the absolute increase in skewness, while $[\underline\eta,\bar\eta]$ bounds this increase relative to the number of added coupling edges.
Underlines and overlines denote lower and upper bounds, respectively.
\Cref{tab:eval-skewness} summarizes the results.

\begin{table}[!htbp]
  \centering
  \scriptsize
  \setlength{\tabcolsep}{1.5pt}
  \caption{Skewness over 100 matched inputs per domain. First five numeric columns: shifted geometric mean bounds; last two: counts with $\bar a=0$ and $L_C=U_C$, respectively. Bold: our method.}
  \label{tab:eval-skewness}
  \begin{tabular}{llrrrrrrr}
  \toprule
  Domain & Method & $U_C$ &
  $\underline{a}$ &
  $\bar{a}$ &
  $\underline{\eta}$ &
  $\bar{\eta}$ &
  \makecell{No added\\skewness} &
  \makecell{Exact\\skewness} \\
  \midrule
PGLib & R & 21.502 & 1.014 & 9.637 & 0.108 & 0.677 & 0 & 14 \\
 & HD & 20.617 & 1.106 & 9.128 & 0.113 & 0.645 & 0 & 16 \\
 & KK & 17.656 & 0.861 & 6.907 & 0.082 & 0.513 & 3 & 21 \\
 & \textbf{OF-H} & 11.377 & 0.057 & 1.928 & 0.005 & 0.213 & 37 & 46 \\
 & \textbf{OF-B} & 10.916 & 0.000 & 1.454 & 0.000 & 0.176 & 48 & 49 \\
\midrule
SNDlib & R & 3.479 & 2.806 & 2.917 & 0.290 & 0.299 & 8 & 90 \\
 & HD & 3.106 & 2.542 & 2.563 & 0.259 & 0.262 & 10 & 96 \\
 & KK & 2.208 & 1.645 & 1.649 & 0.175 & 0.176 & 18 & 99 \\
 & \textbf{OF-H} & 0.584 & 0.014 & 0.014 & 0.001 & 0.001 & 98 & 100 \\
 & \textbf{OF-B} & 0.577 & 0.000 & 0.000 & 0.000 & 0.000 & 100 & 100 \\
\bottomrule
\end{tabular}

\end{table}

\FloatBarrier

\subsection{Other structural characteristics}\label{app:eval-structure}

\begin{figure}[!htbp]
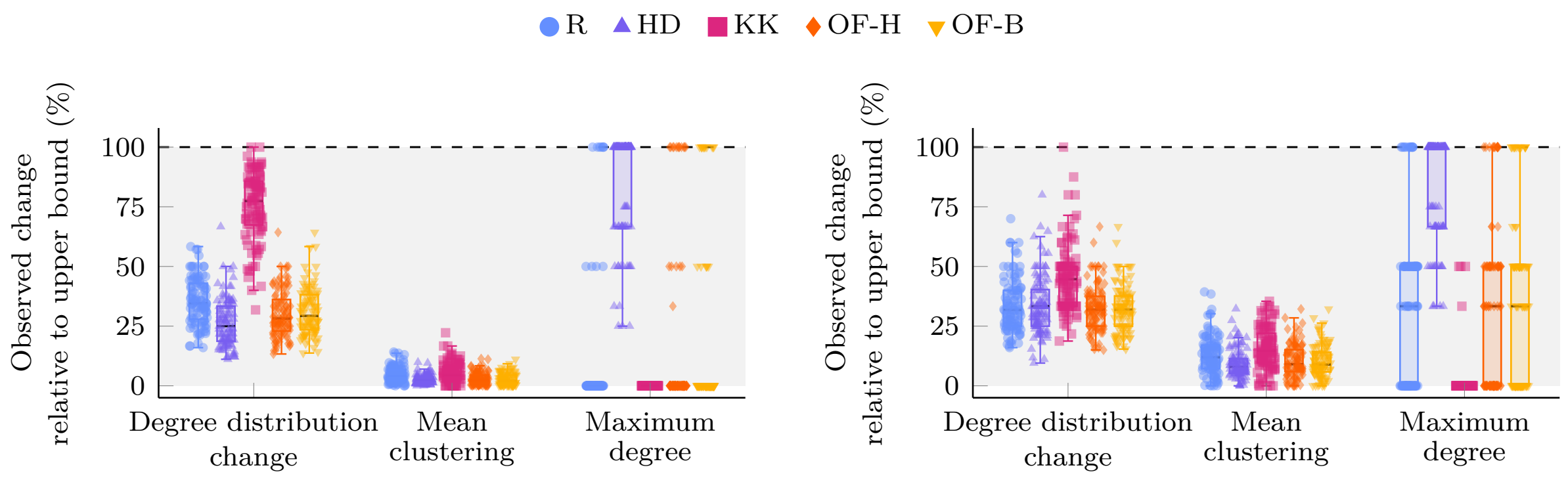

  \centering
  \pgfplotslegendfromname{boundedlegend}\par\smallskip
  \begin{subfigure}[t]{0.485\linewidth}
    \centering
    \input{resources/evaluation/fig-eval-bounded-pglib.tikz}
    \caption{PGLib}
  \end{subfigure}
  \hfill
  \begin{subfigure}[t]{0.485\linewidth}
    \centering
    \input{resources/evaluation/fig-eval-bounded-snd.tikz}
    \caption{SNDlib}
  \end{subfigure}
  \caption{Structural changes that can be bounded (see \Cref{prop:preserved}). Points: composites; boxes: quartiles and median; whiskers: observations within 1.5 interquartile ranges (IQRs); dashed line: upper bound.}
  \label{fig:eval-bounded}
\end{figure}

\begin{figure}[!htbp]
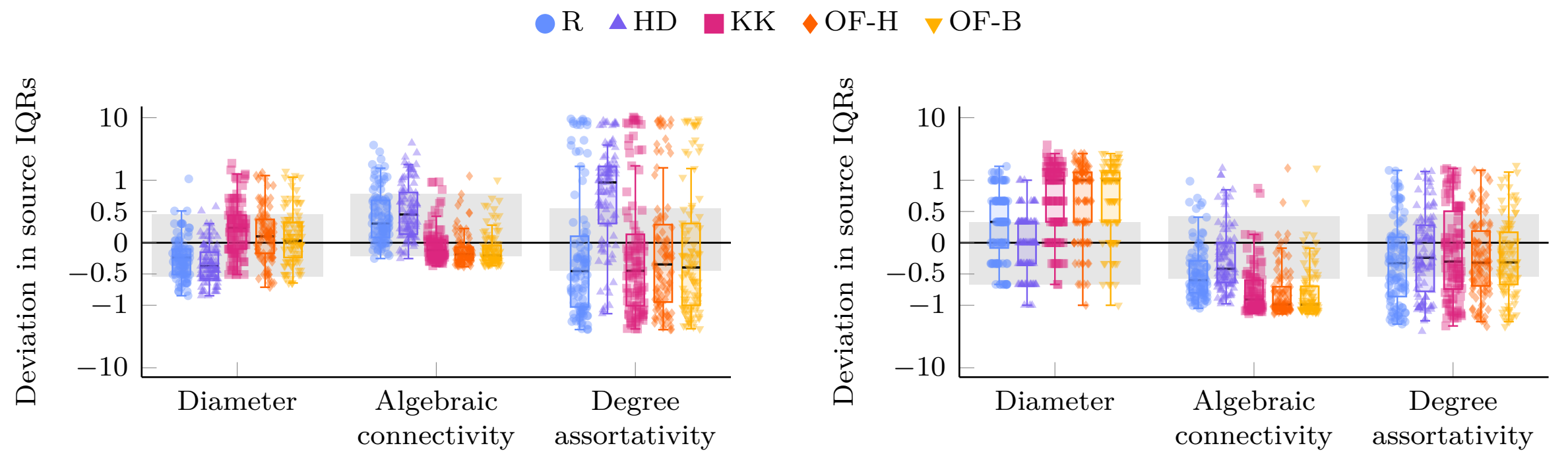

  \centering
  \pgfplotslegendfromname{structurelegend}\par\smallskip
  \begin{subfigure}[t]{0.485\linewidth}
    \centering
    \begingroup
      \input{resources/evaluation/fig-eval-structure-pglib.tikz}
    \endgroup
    \caption{PGLib}
  \end{subfigure}
  \hfill
  \begin{subfigure}[t]{0.485\linewidth}
    \centering
    \begingroup
      \input{resources/evaluation/fig-eval-structure-snd.tikz}
    \endgroup
    \caption{SNDlib}
  \end{subfigure}
  \caption{Structural changes without general bounds. Values centered at source medians and scaled by source IQRs. Shading: source quartiles; symmetric logarithmic axis, linear on $[-1,1]$. Boxplots as in \Cref{fig:eval-bounded}.}
  \label{fig:eval-structure}
\end{figure}

We discuss \Cref{fig:eval-bounded,fig:eval-structure} in \Cref{sec:eval:composite}.

\FloatBarrier

\FloatBarrier

\subsection{Degrees of selected merge vertices}\label{app:merge-vertex-degrees}

\begin{figure}[htbp]
  \centering
\begingroup
\pgfplotslegendfromname{mergedegreelegend}\par\smallskip
\begin{subfigure}[t]{0.485\linewidth}
\centering
\begin{tikzpicture}
\begin{axis}[
  width=.97\linewidth,height=.5567\linewidth,
  axis lines*=left,axis line style={line width=.5pt},
  xmin=.5,xmax=17.5,ymin=-.015,ymax=.78,
  xtick={1,3,5,7,9,11,13,15,17},
  ytick={0,.2,.4,.6},yticklabels={0,20,40,60},
  tick label style={font=\scriptsize},label style={font=\scriptsize},
  ylabel={\shortstack{Selected merge\\vertices [\%]}},
  xlabel={Degree before the join},
  ymajorgrids=true,grid style={gray!20,line width=.3pt},
  legend columns=5,legend cell align=left,
  legend style={font=\scriptsize,draw=none,
    /tikz/every even column/.append style={column sep=3mm}},
  every axis plot/.append style={line width=.8pt,mark size=1.7pt},
  legend to name=mergedegreelegend,
]
\addplot+[color=pigMethod0,mark=o] table[x=degree,y=fraction] {resources/evaluation/plot-data/merge-node-degree-pglib-random.dat};
\addlegendentry{R}
\addplot+[color=pigMethod1,mark=triangle] table[x=degree,y=fraction] {resources/evaluation/plot-data/merge-node-degree-pglib-high_degree.dat};
\addlegendentry{HD}
\addplot+[color=pigMethod2,mark=square] table[x=degree,y=fraction] {resources/evaluation/plot-data/merge-node-degree-pglib-kamada_kawai.dat};
\addlegendentry{KK}
\addplot+[color=pigMethod3,mark=diamond,dashed] table[x=degree,y=fraction] {resources/evaluation/plot-data/merge-node-degree-pglib-heuristic_outer_face.dat};
\addlegendentry{OF-H}
\addplot+[color=pigMethod4,mark=triangle,mark options={rotate=180}] table[x=degree,y=fraction] {resources/evaluation/plot-data/merge-node-degree-pglib-exact_outer_face.dat};
\addlegendentry{OF-B}
\end{axis}
\end{tikzpicture}
\caption{PGLib}
\end{subfigure}
\hfill
\begin{subfigure}[t]{0.485\linewidth}
\centering
\begin{tikzpicture}
\begin{axis}[
  width=.97\linewidth,height=.5567\linewidth,
  axis lines*=left,axis line style={line width=.5pt},
  xmin=.5,xmax=17.5,ymin=-.015,ymax=.78,
  xtick={1,3,5,7,9,11,13,15,17},
  ytick={0,.2,.4,.6},yticklabels={0,20,40,60},
  tick label style={font=\scriptsize},label style={font=\scriptsize},
  ylabel={\shortstack{Selected merge\\vertices [\%]}},
  xlabel={Degree before the join},
  ymajorgrids=true,grid style={gray!20,line width=.3pt},
  legend columns=5,legend cell align=left,
  legend style={font=\scriptsize,draw=none,
    /tikz/every even column/.append style={column sep=3mm}},
  every axis plot/.append style={line width=.8pt,mark size=1.7pt},
]
\addplot+[color=pigMethod0,mark=o,forget plot] table[x=degree,y=fraction] {resources/evaluation/plot-data/merge-node-degree-snd-random.dat};
\addplot+[color=pigMethod1,mark=triangle,forget plot] table[x=degree,y=fraction] {resources/evaluation/plot-data/merge-node-degree-snd-high_degree.dat};
\addplot+[color=pigMethod2,mark=square,forget plot] table[x=degree,y=fraction] {resources/evaluation/plot-data/merge-node-degree-snd-kamada_kawai.dat};
\addplot+[color=pigMethod3,mark=diamond,dashed,forget plot] table[x=degree,y=fraction] {resources/evaluation/plot-data/merge-node-degree-snd-heuristic_outer_face.dat};
\addplot+[color=pigMethod4,mark=triangle,mark options={rotate=180},forget plot] table[x=degree,y=fraction] {resources/evaluation/plot-data/merge-node-degree-snd-exact_outer_face.dat};
\end{axis}
\end{tikzpicture}
\caption{SNDlib}
\end{subfigure}
\endgroup
  \caption{Degrees before each join; 3,218 coupling edge endpoints per method for PGLib and 2,156 for SNDlib.}
  \label{fig:app-merge-vertex-degrees}
\end{figure}
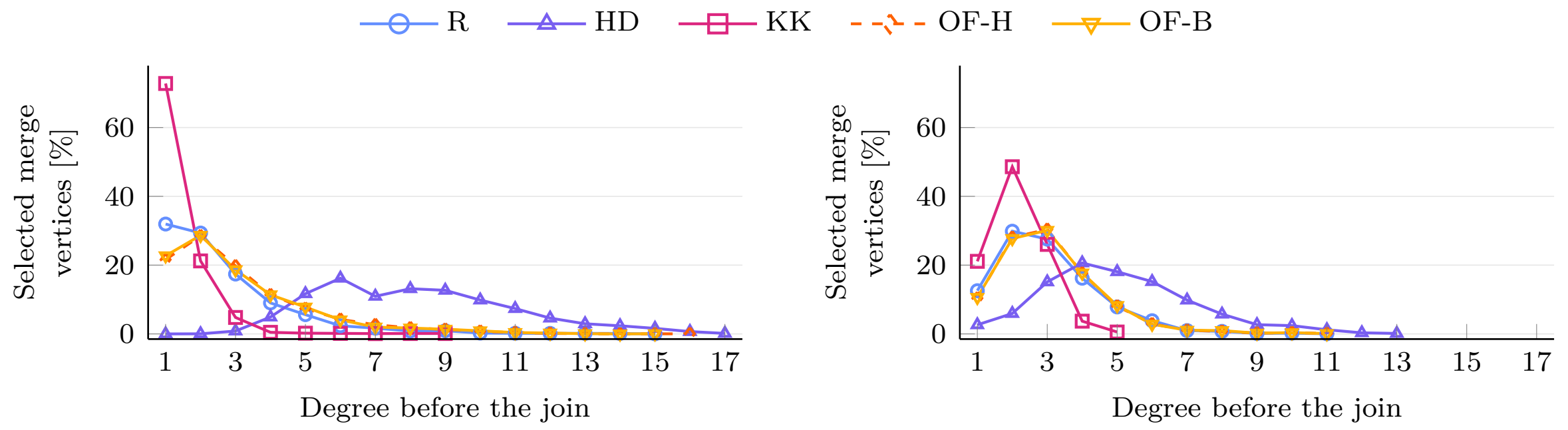

KK most often selects endpoints of degree one in PGLib and degree two in SNDlib (\Cref{fig:app-merge-vertex-degrees}).

\FloatBarrier

\subsection{Generation and application outcomes}\label{app:eval-funnel}

\begin{figure}[htbp]
  \centering
  \begingroup
\newcommand{\funnelflow}[6]{%
  \draw[draw=#6, line width=#5, line cap=butt]
    (#1,#2) .. controls +(0.65,0) and +(-0.65,0) .. (#3,#4);%
}
\newcommand{\funnelnosolution}[5]{%
  \pgfmathsetlengthmacro{\funnelInnerWidth}{#5-.5pt}%
  \draw[draw=black!60,double=black!12,double distance=\funnelInnerWidth,
    line width=.25pt,line cap=butt]
    (#1,#2) .. controls +(0.65,0) and +(-0.65,0) .. (#3,#4);%
}
\newcommand{\funnelexpand}[2]{%
  \draw[draw=black!72,line width=11pt,line cap=butt]
    (2.20,#1) -- (2.55,#1);
  \draw[draw=black!72,line width=30pt,line cap=butt]
    (3.15,#1) -- (#2,#1);
  \path[draw=black,fill=white,line width=.5pt,line join=round]
    ([yshift=5.25pt]2.55,#1)
    .. controls +(0.20,0) and +(-0.20,0) .. ([yshift=14.75pt]3.15,#1)
    -- ([yshift=-14.75pt]3.15,#1)
    .. controls +(-0.20,0) and +(0.20,0) .. ([yshift=-5.25pt]2.55,#1)
    -- cycle;
  \node at (2.85,#1) {\(\times 5\)};
}
\begin{subfigure}[t]{.49\linewidth}
  \centering
  \begin{tikzpicture}[font=\fontsize{7.5}{9}\selectfont, every node/.style={inner sep=1pt}]
\path[use as bounding box] (-.10,7.25) rectangle (7.60,12.10);
\node[anchor=south,align=center] at (.65,11.00) {262 attempts\\(100\%)};
\funnelflow{.60}{10.45}{2.20}{10.45}{11pt}{black!72}
\node[anchor=south,align=center] at (2.60,11.00) {100 accepted\\(38.2\%)};
\funnelflow{.60}{10.193}{2.20}{9.05}{3.63pt}{black!18}
\node[anchor=west,align=left] at (2.35,9.05) {Feasibility check\\33 (12.6\%)};
\funnelflow{.60}{9.880}{2.20}{7.90}{14.19pt}{black!50}
\node[anchor=west,align=left] at (2.35,7.90) {Size requirements\\unmet\\129 (49.2\%)};
\node[anchor=west] at (2.35,9.70) {Rejected};
\funnelexpand{10.45}{4.30}
\node[font=\fontsize{7.5}{9}\selectfont\bfseries] at (5.85,11.75) {SCOTS, \(\TAL=1\) (\(N=500\))};
\funnelflow{4.30}{10.668}{5.60}{11.00}{17.58pt}{black!82}
\node[anchor=west,align=left,text width=1.85cm] at (5.75,11.00) {Optimal\\293 (58.6\%)};
\funnelflow{4.30}{10.322}{5.60}{10.00}{2.1pt}{black!58}
\node[anchor=west,align=left,text width=1.85cm] at (5.75,10.00) {Feasible\\35 (7.0\%)};
\funnelnosolution{4.30}{10.267}{5.60}{9.00}{1.02pt}
\node[anchor=west,align=left,text width=1.85cm] at (5.75,9.00) {No solution\\17 (3.4\%)};
\funnelflow{4.30}{10.086}{5.60}{8.00}{9.3pt}{black!34}
\node[anchor=west,align=left,text width=1.85cm] at (5.75,8.00) {Infeasible\\155 (31.0\%)};
\node[align=center,text=white] at (3.725,10.45) {500\\(100\%)};
\end{tikzpicture}
  \caption{PGLib}
\end{subfigure}\hfill
\begin{subfigure}[t]{.49\linewidth}
  \centering
  \begin{tikzpicture}[font=\fontsize{7.5}{9}\selectfont, every node/.style={inner sep=1pt}]
\path[use as bounding box] (-.10,7.25) rectangle (7.60,12.10);
\node[anchor=south,align=center] at (.65,11.00) {142 attempts\\(100\%)};
\funnelflow{.60}{10.45}{2.20}{10.45}{11pt}{black!72}
\node[anchor=south,align=center] at (2.60,11.00) {100 accepted\\(70.4\%)};
\funnelflow{.60}{10.176}{2.20}{7.90}{4.62pt}{black!50}
\node[anchor=west,align=left] at (2.35,7.90) {Size requirements\\unmet\\42 (29.6\%)};
\node[anchor=west] at (2.35,8.70) {Rejected};
\funnelexpand{10.45}{4.30}
\node[font=\fontsize{7.5}{9}\selectfont\bfseries] at (5.85,11.75) {SND (\(N=500\))};
\funnelflow{4.30}{10.840}{5.60}{11.00}{7.8pt}{black!82}
\node[anchor=west,align=left,text width=1.85cm] at (5.75,11.00) {Optimal\\130 (26.0\%)};
\funnelflow{4.30}{10.614}{5.60}{10.00}{5.04pt}{black!58}
\node[anchor=west,align=left,text width=1.85cm] at (5.75,10.00) {Feasible\\84 (16.8\%)};
\funnelnosolution{4.30}{10.224}{5.60}{9.00}{17.16pt}
\node[anchor=west,align=left,text width=1.85cm] at (5.75,9.00) {No solution\\286 (57.2\%)};
\node[anchor=west,align=left,text width=1.85cm] at (5.75,8.00) {Infeasible\\0 (0.0\%)};
\node[align=center,text=white] at (3.725,10.45) {500\\(100\%)};
\end{tikzpicture}
  \caption{SNDlib}
\end{subfigure}
\endgroup
  \caption{Percentages refer to all attempts or 500 composites per application. Size requirements concern instance size and available outer face vertices. Band widths are proportional within stages.}
  \label{fig:app-evaluation-funnel}
\end{figure}
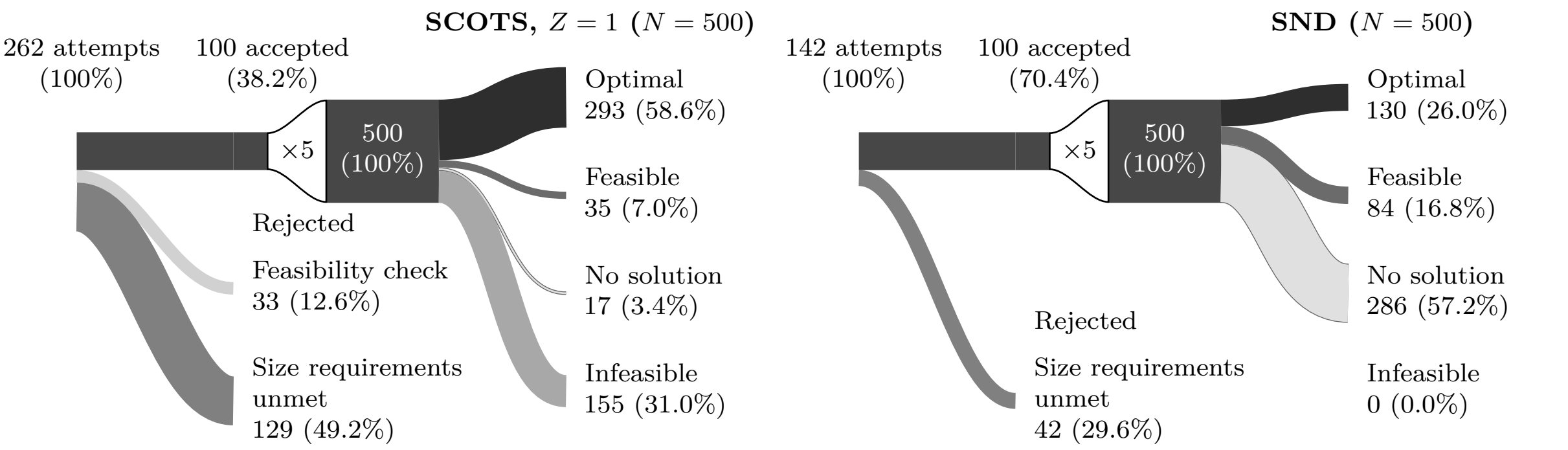

In \Cref{fig:app-evaluation-funnel}, ``Feasible'' denotes a run with a feasible solution but no proof of optimality.
``No solution'' denotes a run ending without either a feasible solution or a proof of infeasibility.

The PGLib rejection counts depend on the fixed evaluation order R, HD, KK, OF-H, OF-B, since each attempt stops at the first failure.
All SNDlib rejections occur at the target size check.

\begin{figure}[!htbp]
  \centering
\begingroup
\tikzset{
  pig infeasible/.style={fill=white,postaction={pattern=dots,pattern color=#1}},
  pig time limit/.style={fill=white,postaction={pattern=north east lines,pattern color=#1}},
}
\pgfplotslegendfromname{applicationoutcomelegend}\par\smallskip
\begin{subfigure}[t]{0.485\linewidth}
  \centering
\begingroup
\pigSetRowSpacing
\begin{tikzpicture}
\begin{axis}[
  set layers=standard,clip=false,xbar stacked,/pgf/bar width=\pigBarWidth,
  ymin=\pigSourceYMin,ymax=4.64,ytick={\pigSourceY,0,1,2,3,4},y dir=reverse,
  yticklabels={Source,R,HD,KK,OF-H,OF-B},
  width=.97\linewidth,height=\pigSourcePlotHeight,
  xmin=0,xmax=100,xtick={0,25,50,75,100},
  axis lines*=left,axis line style={line width=.5pt},
  tick label style={font=\scriptsize},label style={font=\scriptsize},
  xlabel={Number of instances},xmajorgrids=true,grid style={gray!18,line width=.3pt},
  legend style={font=\scriptsize,legend columns=3,draw=none,
    /tikz/every even column/.append style={column sep=3mm}},
  legend to name=applicationoutcomelegend,
]
\addlegendimage{area legend,draw=black,fill=black}
\addlegendentry{Feasible}
\addlegendimage{area legend,draw=black,pig infeasible=black}
\addlegendentry{Infeasible}
\addlegendimage{area legend,draw=black,fill=white}
\addlegendentry{No solution}
\addplot[draw=pigMethod5,fill=pigMethod5,forget plot] coordinates {(12,\pigSourceY) (0,0) (0,1) (0,2) (0,3) (0,4)};
\addplot[draw=pigMethod5,pig infeasible=pigMethod5,forget plot] coordinates {(4,\pigSourceY) (0,0) (0,1) (0,2) (0,3) (0,4)};
\addplot[draw=pigMethod5,fill=white,forget plot] coordinates {(0,\pigSourceY) (0,0) (0,1) (0,2) (0,3) (0,4)};
\addplot[draw=pigMethod0,fill=pigMethod0,forget plot] coordinates {(0,\pigSourceY) (69,0) (0,1) (0,2) (0,3) (0,4)};
\addplot[draw=pigMethod0,pig infeasible=pigMethod0,forget plot] coordinates {(0,\pigSourceY) (27,0) (0,1) (0,2) (0,3) (0,4)};
\addplot[draw=pigMethod0,fill=white,forget plot] coordinates {(0,\pigSourceY) (4,0) (0,1) (0,2) (0,3) (0,4)};
\addplot[draw=pigMethod1,fill=pigMethod1,forget plot] coordinates {(0,\pigSourceY) (0,0) (70,1) (0,2) (0,3) (0,4)};
\addplot[draw=pigMethod1,pig infeasible=pigMethod1,forget plot] coordinates {(0,\pigSourceY) (0,0) (25,1) (0,2) (0,3) (0,4)};
\addplot[draw=pigMethod1,fill=white,forget plot] coordinates {(0,\pigSourceY) (0,0) (5,1) (0,2) (0,3) (0,4)};
\addplot[draw=pigMethod2,fill=pigMethod2,forget plot] coordinates {(0,\pigSourceY) (0,0) (0,1) (54,2) (0,3) (0,4)};
\addplot[draw=pigMethod2,pig infeasible=pigMethod2,forget plot] coordinates {(0,\pigSourceY) (0,0) (0,1) (43,2) (0,3) (0,4)};
\addplot[draw=pigMethod2,fill=white,forget plot] coordinates {(0,\pigSourceY) (0,0) (0,1) (3,2) (0,3) (0,4)};
\addplot[draw=pigMethod3,fill=pigMethod3,forget plot] coordinates {(0,\pigSourceY) (0,0) (0,1) (0,2) (67,3) (0,4)};
\addplot[draw=pigMethod3,pig infeasible=pigMethod3,forget plot] coordinates {(0,\pigSourceY) (0,0) (0,1) (0,2) (30,3) (0,4)};
\addplot[draw=pigMethod3,fill=white,forget plot] coordinates {(0,\pigSourceY) (0,0) (0,1) (0,2) (3,3) (0,4)};
\addplot[draw=pigMethod4,fill=pigMethod4,forget plot] coordinates {(0,\pigSourceY) (0,0) (0,1) (0,2) (0,3) (68,4)};
\addplot[draw=pigMethod4,pig infeasible=pigMethod4,forget plot] coordinates {(0,\pigSourceY) (0,0) (0,1) (0,2) (0,3) (30,4)};
\addplot[draw=pigMethod4,fill=white,forget plot] coordinates {(0,\pigSourceY) (0,0) (0,1) (0,2) (0,3) (2,4)};
\pgfplotsextra{
\begin{pgfonlayer}{axis foreground}
\node[font=\scriptsize,text=white,inner sep=0pt] at (axis cs:6,\pigSourceY) {12};
\node[font=\scriptsize,text=black,fill=white,inner sep=1pt] at (axis cs:14,\pigSourceY) {4};
\node[font=\scriptsize,text=white,inner sep=0pt] at (axis cs:34.5,0) {69};
\node[font=\scriptsize,text=black,fill=white,inner sep=1pt] at (axis cs:82.5,0) {27};
\node[font=\scriptsize,text=black,fill=white,inner sep=1pt] at (axis cs:98,0) {4};
\node[font=\scriptsize,text=white,inner sep=0pt] at (axis cs:35,1) {70};
\node[font=\scriptsize,text=black,fill=white,inner sep=1pt] at (axis cs:82.5,1) {25};
\node[font=\scriptsize,text=black,fill=white,inner sep=1pt] at (axis cs:97.5,1) {5};
\node[font=\scriptsize,text=white,inner sep=0pt] at (axis cs:27,2) {54};
\node[font=\scriptsize,text=black,fill=white,inner sep=1pt] at (axis cs:75.5,2) {43};
\node[font=\scriptsize,text=black,fill=white,inner sep=1pt] at (axis cs:98.5,2) {3};
\node[font=\scriptsize,text=white,inner sep=0pt] at (axis cs:33.5,3) {67};
\node[font=\scriptsize,text=black,fill=white,inner sep=1pt] at (axis cs:82,3) {30};
\node[font=\scriptsize,text=black,fill=white,inner sep=1pt] at (axis cs:98.5,3) {3};
\node[font=\scriptsize,text=white,inner sep=0pt] at (axis cs:34,4) {68};
\node[font=\scriptsize,text=black,fill=white,inner sep=1pt] at (axis cs:83,4) {30};
\node[font=\scriptsize,text=black,fill=white,inner sep=1pt] at (axis cs:99,4) {2};
\end{pgfonlayer}
}
\end{axis}
\end{tikzpicture}
\endgroup
  \caption{SCOTS}
\end{subfigure}
\hfill
\begin{subfigure}[t]{0.485\linewidth}
  \centering
\begingroup
\pigSetRowSpacing
\begin{tikzpicture}
\begin{axis}[
  set layers=standard,clip=false,xbar stacked,/pgf/bar width=\pigBarWidth,
  ymin=\pigSourceYMin,ymax=4.64,ytick={\pigSourceY,0,1,2,3,4},y dir=reverse,
  yticklabels={Source,R,HD,KK,OF-H,OF-B},
  width=.97\linewidth,height=\pigSourcePlotHeight,
  xmin=0,xmax=100,xtick={0,25,50,75,100},
  axis lines*=left,axis line style={line width=.5pt},
  tick label style={font=\scriptsize},label style={font=\scriptsize},
  xlabel={Number of instances},xmajorgrids=true,grid style={gray!18,line width=.3pt},
  legend style={font=\scriptsize,legend columns=3,draw=none,
    /tikz/every even column/.append style={column sep=3mm}},
]
\addplot[draw=pigMethod5,fill=pigMethod5,forget plot] coordinates {(14,\pigSourceY) (0,0) (0,1) (0,2) (0,3) (0,4)};
\addplot[draw=pigMethod5,pig infeasible=pigMethod5,forget plot] coordinates {(0,\pigSourceY) (0,0) (0,1) (0,2) (0,3) (0,4)};
\addplot[draw=pigMethod5,fill=white,forget plot] coordinates {(2,\pigSourceY) (0,0) (0,1) (0,2) (0,3) (0,4)};
\addplot[draw=pigMethod0,fill=pigMethod0,forget plot] coordinates {(0,\pigSourceY) (40,0) (0,1) (0,2) (0,3) (0,4)};
\addplot[draw=pigMethod0,pig infeasible=pigMethod0,forget plot] coordinates {(0,\pigSourceY) (0,0) (0,1) (0,2) (0,3) (0,4)};
\addplot[draw=pigMethod0,fill=white,forget plot] coordinates {(0,\pigSourceY) (60,0) (0,1) (0,2) (0,3) (0,4)};
\addplot[draw=pigMethod1,fill=pigMethod1,forget plot] coordinates {(0,\pigSourceY) (0,0) (47,1) (0,2) (0,3) (0,4)};
\addplot[draw=pigMethod1,pig infeasible=pigMethod1,forget plot] coordinates {(0,\pigSourceY) (0,0) (0,1) (0,2) (0,3) (0,4)};
\addplot[draw=pigMethod1,fill=white,forget plot] coordinates {(0,\pigSourceY) (0,0) (53,1) (0,2) (0,3) (0,4)};
\addplot[draw=pigMethod2,fill=pigMethod2,forget plot] coordinates {(0,\pigSourceY) (0,0) (0,1) (39,2) (0,3) (0,4)};
\addplot[draw=pigMethod2,pig infeasible=pigMethod2,forget plot] coordinates {(0,\pigSourceY) (0,0) (0,1) (0,2) (0,3) (0,4)};
\addplot[draw=pigMethod2,fill=white,forget plot] coordinates {(0,\pigSourceY) (0,0) (0,1) (61,2) (0,3) (0,4)};
\addplot[draw=pigMethod3,fill=pigMethod3,forget plot] coordinates {(0,\pigSourceY) (0,0) (0,1) (0,2) (44,3) (0,4)};
\addplot[draw=pigMethod3,pig infeasible=pigMethod3,forget plot] coordinates {(0,\pigSourceY) (0,0) (0,1) (0,2) (0,3) (0,4)};
\addplot[draw=pigMethod3,fill=white,forget plot] coordinates {(0,\pigSourceY) (0,0) (0,1) (0,2) (56,3) (0,4)};
\addplot[draw=pigMethod4,fill=pigMethod4,forget plot] coordinates {(0,\pigSourceY) (0,0) (0,1) (0,2) (0,3) (44,4)};
\addplot[draw=pigMethod4,pig infeasible=pigMethod4,forget plot] coordinates {(0,\pigSourceY) (0,0) (0,1) (0,2) (0,3) (0,4)};
\addplot[draw=pigMethod4,fill=white,forget plot] coordinates {(0,\pigSourceY) (0,0) (0,1) (0,2) (0,3) (56,4)};
\pgfplotsextra{
\begin{pgfonlayer}{axis foreground}
\node[font=\scriptsize,text=white,inner sep=0pt] at (axis cs:7,\pigSourceY) {14};
\node[font=\scriptsize,text=black,fill=white,inner sep=1pt,anchor=south,yshift=6pt] at (axis cs:15,\pigSourceY) {2};
\node[font=\scriptsize,text=white,inner sep=0pt] at (axis cs:20,0) {40};
\node[font=\scriptsize,text=black,fill=white,inner sep=1pt] at (axis cs:70,0) {60};
\node[font=\scriptsize,text=white,inner sep=0pt] at (axis cs:23.5,1) {47};
\node[font=\scriptsize,text=black,fill=white,inner sep=1pt] at (axis cs:73.5,1) {53};
\node[font=\scriptsize,text=white,inner sep=0pt] at (axis cs:19.5,2) {39};
\node[font=\scriptsize,text=black,fill=white,inner sep=1pt] at (axis cs:69.5,2) {61};
\node[font=\scriptsize,text=white,inner sep=0pt] at (axis cs:22,3) {44};
\node[font=\scriptsize,text=black,fill=white,inner sep=1pt] at (axis cs:72,3) {56};
\node[font=\scriptsize,text=white,inner sep=0pt] at (axis cs:22,4) {44};
\node[font=\scriptsize,text=black,fill=white,inner sep=1pt] at (axis cs:72,4) {56};
\end{pgfonlayer}
}
\end{axis}
\end{tikzpicture}
\endgroup
  \caption{SND}
\end{subfigure}
\endgroup
  \caption{Application outcomes; 100 composites per method and 16 sources per application; \SI{10800}{\second} solver limit. ``Feasible'' includes optimal solutions. ``Infeasible'' requires a solver proof; ``No solution'' denotes runs ending without either a feasible solution or an infeasibility proof.}
  \label{fig:eval-application-outcomes}
\end{figure}
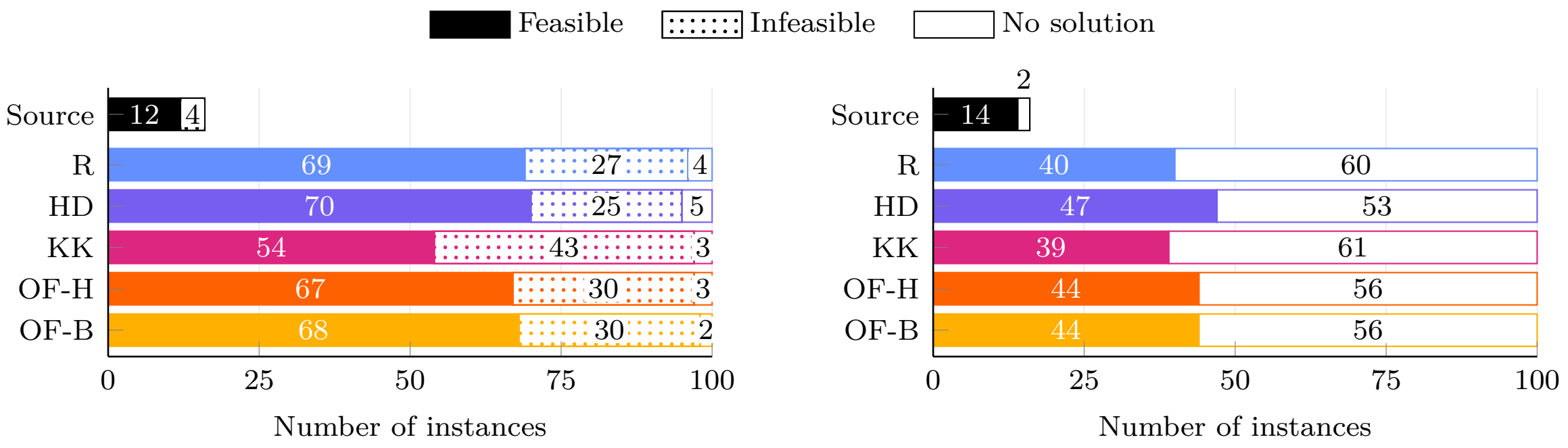

\FloatBarrier
\subsection{Flows between subinstances}\label{app:eval-flows}

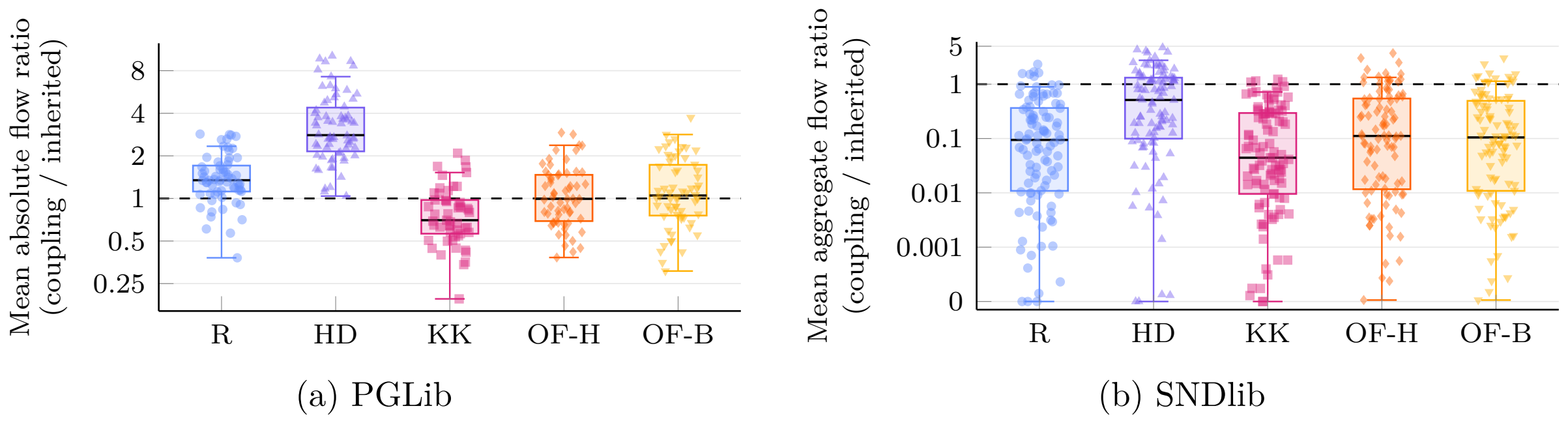
\begin{figure}[!htbp]
  \centering
  \begin{subfigure}[b]{0.485\linewidth}
    \centering
    \makebox[\linewidth][c]{
\begin{tikzpicture}
\begin{axis}[
  width=.97\linewidth,height=.5567\linewidth,
  xmin=-.55,xmax=4.55,ymin=.16,ymax=12.5,
  ymode=log,log basis y=2,
  xtick={0,1,2,3,4},
  xticklabels={R,HD,KK,OF-H,OF-B},
  ytick={.25,.5,1,2,4,8},yticklabels={$0.25$,$0.5$,$1$,$2$,$4$,$8$},
  minor y tick num=0,
  axis lines*=left,axis line style={line width=.5pt},
  tick label style={font=\scriptsize},label style={font=\scriptsize},
  ylabel={\shortstack{Mean absolute flow ratio\\(coupling / inherited)}},
  ymajorgrids=true,grid style={gray!20,line width=.3pt},
]
\draw[black,dashed,line width=.55pt] (axis cs:-.55,1) -- (axis cs:4.55,1);
\addplot[only marks,mark=*,mark size=1.15pt,color=pigMethod0,opacity=.45] table[x=x,y=y] {resources/evaluation/plot-data/scopf-flow-random.dat};
\draw[color=pigMethod0,line width=.45pt] (axis cs:0,0.380500492541) -- (axis cs:0,1.11903477163);
\draw[color=pigMethod0,line width=.45pt] (axis cs:0,1.707860891) -- (axis cs:0,2.33843435239);
\draw[color=pigMethod0,line width=.45pt] (axis cs:-0.13,0.380500492541) -- (axis cs:0.13,0.380500492541);
\draw[color=pigMethod0,line width=.45pt] (axis cs:-0.13,2.33843435239) -- (axis cs:0.13,2.33843435239);
\path[draw=pigMethod0,fill=pigMethod0,fill opacity=.16,draw opacity=1,line width=.5pt] (axis cs:-0.25,1.11903477163) rectangle (axis cs:0.25,1.707860891);
\draw[black,line width=.65pt] (axis cs:-0.25,1.34508860579) -- (axis cs:0.25,1.34508860579);
\addplot[only marks,mark=triangle*,mark size=1.15pt,color=pigMethod1,opacity=.45] table[x=x,y=y] {resources/evaluation/plot-data/scopf-flow-high_degree.dat};
\draw[color=pigMethod1,line width=.45pt] (axis cs:1,1.03868481317) -- (axis cs:1,2.15193785895);
\draw[color=pigMethod1,line width=.45pt] (axis cs:1,4.38808256049) -- (axis cs:1,7.26212403521);
\draw[color=pigMethod1,line width=.45pt] (axis cs:0.87,1.03868481317) -- (axis cs:1.13,1.03868481317);
\draw[color=pigMethod1,line width=.45pt] (axis cs:0.87,7.26212403521) -- (axis cs:1.13,7.26212403521);
\path[draw=pigMethod1,fill=pigMethod1,fill opacity=.16,draw opacity=1,line width=.5pt] (axis cs:0.75,2.15193785895) rectangle (axis cs:1.25,4.38808256049);
\draw[black,line width=.65pt] (axis cs:0.75,2.80139287931) -- (axis cs:1.25,2.80139287931);
\addplot[only marks,mark=square*,mark size=1.15pt,color=pigMethod2,opacity=.45] table[x=x,y=y] {resources/evaluation/plot-data/scopf-flow-kamada_kawai.dat};
\draw[color=pigMethod2,line width=.45pt] (axis cs:2,0.194975663159) -- (axis cs:2,0.563295547949);
\draw[color=pigMethod2,line width=.45pt] (axis cs:2,0.975107633401) -- (axis cs:2,1.52462294188);
\draw[color=pigMethod2,line width=.45pt] (axis cs:1.87,0.194975663159) -- (axis cs:2.13,0.194975663159);
\draw[color=pigMethod2,line width=.45pt] (axis cs:1.87,1.52462294188) -- (axis cs:2.13,1.52462294188);
\path[draw=pigMethod2,fill=pigMethod2,fill opacity=.16,draw opacity=1,line width=.5pt] (axis cs:1.75,0.563295547949) rectangle (axis cs:2.25,0.975107633401);
\draw[black,line width=.65pt] (axis cs:1.75,0.701679143976) -- (axis cs:2.25,0.701679143976);
\addplot[only marks,mark=diamond*,mark size=1.15pt,color=pigMethod3,opacity=.45] table[x=x,y=y] {resources/evaluation/plot-data/scopf-flow-heuristic_outer_face.dat};
\draw[color=pigMethod3,line width=.45pt] (axis cs:3,0.38240753699) -- (axis cs:3,0.691863466199);
\draw[color=pigMethod3,line width=.45pt] (axis cs:3,1.46861571535) -- (axis cs:3,2.37581420107);
\draw[color=pigMethod3,line width=.45pt] (axis cs:2.87,0.38240753699) -- (axis cs:3.13,0.38240753699);
\draw[color=pigMethod3,line width=.45pt] (axis cs:2.87,2.37581420107) -- (axis cs:3.13,2.37581420107);
\path[draw=pigMethod3,fill=pigMethod3,fill opacity=.16,draw opacity=1,line width=.5pt] (axis cs:2.75,0.691863466199) rectangle (axis cs:3.25,1.46861571535);
\draw[black,line width=.65pt] (axis cs:2.75,0.99353072189) -- (axis cs:3.25,0.99353072189);
\addplot[only marks,mark=triangle*,mark size=1.15pt,color=pigMethod4,opacity=.45,mark options={rotate=180}] table[x=x,y=y] {resources/evaluation/plot-data/scopf-flow-exact_outer_face.dat};
\draw[color=pigMethod4,line width=.45pt] (axis cs:4,0.306714861629) -- (axis cs:4,0.756847667624);
\draw[color=pigMethod4,line width=.45pt] (axis cs:4,1.72859772073) -- (axis cs:4,2.83079815798);
\draw[color=pigMethod4,line width=.45pt] (axis cs:3.87,0.306714861629) -- (axis cs:4.13,0.306714861629);
\draw[color=pigMethod4,line width=.45pt] (axis cs:3.87,2.83079815798) -- (axis cs:4.13,2.83079815798);
\path[draw=pigMethod4,fill=pigMethod4,fill opacity=.16,draw opacity=1,line width=.5pt] (axis cs:3.75,0.756847667624) rectangle (axis cs:4.25,1.72859772073);
\draw[black,line width=.65pt] (axis cs:3.75,1.04890824713) -- (axis cs:4.25,1.04890824713);
\end{axis}
\end{tikzpicture}}
    \caption{PGLib}
    \label{fig:eval-flows-pglib}
  \end{subfigure}
  \hfill
  \begin{subfigure}[b]{0.485\linewidth}
    \centering
    \makebox[\linewidth][c]{
\begin{tikzpicture}
\begin{axis}[
  width=.97\linewidth,height=.5567\linewidth,
  xmin=-.55,xmax=4.55,ymin=-.15,ymax=4.77815125038,
  xtick={0,1,2,3,4},xticklabels={R,HD,KK,OF-H,OF-B},
  ytick={0,1,2,3,4,4.69897000434},yticklabels={$0$,$0.001$,$0.01$,$0.1$,$1$,$5$},
  axis lines*=left,axis line style={line width=.5pt},
  tick label style={font=\scriptsize},label style={font=\scriptsize},
  ylabel={\shortstack{Mean aggregate flow ratio\\(coupling / inherited)}},
  ymajorgrids=true,grid style={gray!20,line width=.3pt},
]
\draw[black,dashed,line width=.55pt] (axis cs:-.55,4) -- (axis cs:4.55,4);
\addplot[only marks,mark=*,mark size=1.15pt,color=pigMethod0,opacity=.45] table[x=x,y=y] {resources/evaluation/plot-data/snd-flow-random.dat};
\draw[color=pigMethod0,line width=.45pt] (axis cs:0,0) -- (axis cs:0,2.03568580918);
\draw[color=pigMethod0,line width=.45pt] (axis cs:0,3.56183690959) -- (axis cs:0,3.95073172165);
\draw[color=pigMethod0,line width=.45pt] (axis cs:-0.13,0) -- (axis cs:0.13,0);
\draw[color=pigMethod0,line width=.45pt] (axis cs:-0.13,3.95073172165) -- (axis cs:0.13,3.95073172165);
\path[draw=pigMethod0,fill=pigMethod0,fill opacity=.16,draw opacity=1,line width=.5pt] (axis cs:-0.25,2.03568580918) rectangle (axis cs:0.25,3.56183690959);
\draw[black,line width=.65pt] (axis cs:-0.25,2.97378676187) -- (axis cs:0.25,2.97378676187);
\addplot[only marks,mark=triangle*,mark size=1.15pt,color=pigMethod1,opacity=.45] table[x=x,y=y] {resources/evaluation/plot-data/snd-flow-high_degree.dat};
\draw[color=pigMethod1,line width=.45pt] (axis cs:1,0) -- (axis cs:1,2.99703184876);
\draw[color=pigMethod1,line width=.45pt] (axis cs:1,4.12021871574) -- (axis cs:1,4.43981505247);
\draw[color=pigMethod1,line width=.45pt] (axis cs:0.87,0) -- (axis cs:1.13,0);
\draw[color=pigMethod1,line width=.45pt] (axis cs:0.87,4.43981505247) -- (axis cs:1.13,4.43981505247);
\path[draw=pigMethod1,fill=pigMethod1,fill opacity=.16,draw opacity=1,line width=.5pt] (axis cs:0.75,2.99703184876) rectangle (axis cs:1.25,4.12021871574);
\draw[black,line width=.65pt] (axis cs:0.75,3.70920722763) -- (axis cs:1.25,3.70920722763);
\addplot[only marks,mark=square*,mark size=1.15pt,color=pigMethod2,opacity=.45] table[x=x,y=y] {resources/evaluation/plot-data/snd-flow-kamada_kawai.dat};
\draw[color=pigMethod2,line width=.45pt] (axis cs:2,0) -- (axis cs:2,1.98051363443);
\draw[color=pigMethod2,line width=.45pt] (axis cs:2,3.46943581339) -- (axis cs:2,3.855516216);
\draw[color=pigMethod2,line width=.45pt] (axis cs:1.87,0) -- (axis cs:2.13,0);
\draw[color=pigMethod2,line width=.45pt] (axis cs:1.87,3.855516216) -- (axis cs:2.13,3.855516216);
\path[draw=pigMethod2,fill=pigMethod2,fill opacity=.16,draw opacity=1,line width=.5pt] (axis cs:1.75,1.98051363443) rectangle (axis cs:2.25,3.46943581339);
\draw[black,line width=.65pt] (axis cs:1.75,2.64632637588) -- (axis cs:2.25,2.64632637588);
\addplot[only marks,mark=diamond*,mark size=1.15pt,color=pigMethod3,opacity=.45] table[x=x,y=y] {resources/evaluation/plot-data/snd-flow-heuristic_outer_face.dat};
\draw[color=pigMethod3,line width=.45pt] (axis cs:3,0.0266445380187) -- (axis cs:3,2.06658888131);
\draw[color=pigMethod3,line width=.45pt] (axis cs:3,3.73647872095) -- (axis cs:3,4.12452485842);
\draw[color=pigMethod3,line width=.45pt] (axis cs:2.87,0.0266445380187) -- (axis cs:3.13,0.0266445380187);
\draw[color=pigMethod3,line width=.45pt] (axis cs:2.87,4.12452485842) -- (axis cs:3.13,4.12452485842);
\path[draw=pigMethod3,fill=pigMethod3,fill opacity=.16,draw opacity=1,line width=.5pt] (axis cs:2.75,2.06658888131) rectangle (axis cs:3.25,3.73647872095);
\draw[black,line width=.65pt] (axis cs:2.75,3.04661557747) -- (axis cs:3.25,3.04661557747);
\addplot[only marks,mark=triangle*,mark size=1.15pt,color=pigMethod4,opacity=.45,mark options={rotate=180}] table[x=x,y=y] {resources/evaluation/plot-data/snd-flow-exact_outer_face.dat};
\draw[color=pigMethod4,line width=.45pt] (axis cs:4,0.0266445380187) -- (axis cs:4,2.03592157202);
\draw[color=pigMethod4,line width=.45pt] (axis cs:4,3.69365883007) -- (axis cs:4,4.05426925563);
\draw[color=pigMethod4,line width=.45pt] (axis cs:3.87,0.0266445380187) -- (axis cs:4.13,0.0266445380187);
\draw[color=pigMethod4,line width=.45pt] (axis cs:3.87,4.05426925563) -- (axis cs:4.13,4.05426925563);
\path[draw=pigMethod4,fill=pigMethod4,fill opacity=.16,draw opacity=1,line width=.5pt] (axis cs:3.75,2.03592157202) rectangle (axis cs:4.25,3.69365883007);
\draw[black,line width=.65pt] (axis cs:3.75,3.02014598179) -- (axis cs:4.25,3.02014598179);
\end{axis}
\end{tikzpicture}}
    \caption{SNDlib}
    \label{fig:eval-flows-snd}
  \end{subfigure}
  \caption{Base-case flows. Boxplots as in \Cref{fig:eval-bounded}. SNDlib axis linear below $10^{-3}$ and logarithmic above.}
  \label{fig:eval-flows}
\end{figure}

For PGLib, we compare the mean absolute base-case flows on coupling and inherited edges in the SCOPF solutions (\Cref{fig:eval-flows-pglib}).

For SNDlib, we solve $\Tract$ with capacities on inherited edges fixed as in preprocessing and one capacity module installed on each coupling edge.
Among minimum-cost routings, we minimize total flow to avoid unnecessary circulation.
We then sum flows over all commodities and both directions on each edge and compare the means on coupling and inherited edges (\Cref{fig:eval-flows-snd}).

\FloatBarrier

\subsection{Outer face vertex selection rules}\label{app:eval-spacing}

For OF-B on PGLib, we compare the four outer face selection rules described in \Cref{sec:eval:sensitivity}, using the main experiment's coupling count and planarization limits.
We solve SCOPF with overload slack to optimality or proven infeasibility (\Cref{tab:app-scopf-sampling}).
The 69 matched inputs are those for which all four rules accept identical subinstances with the same target sizes and seeds.

\begin{table}[htbp]
  \centering
  \scriptsize
  \caption{Counts of SCOPF runs returning a feasible solution with OF-B. Bold: setting used in the main experiment.}
  \label{tab:app-scopf-sampling}
  \setlength{\tabcolsep}{4pt}
  \begin{tabular}{lrr}
  \toprule
  Selection rule & All 100 inputs & 69 matched inputs \\
  \midrule
  \textbf{Consecutive} & 62 & 45 \\
  Equidistant & 60 & 43 \\
  Skip one where possible & 57 & 45 \\
  Consecutive, prefer degree $>1$ & 61 & 46 \\
  \bottomrule
\end{tabular}

\end{table}

\FloatBarrier

\Needspace{0.4\textheight}
\subsection{Normalization across voltage levels}\label{app:eval-normalization}

For PGLib, we normalize the accepted composites before SCOPF evaluation (\Cref{tab:eval-normalization}).
We assume positive total apparent demand in every subinstance.
For each subinstance $i$ with vertex set $V_i$, we use
\[
 \begin{aligned}
   d_i&=\frac{1}{|V_i|}\sum_{v\in V_i}
       \frac{\sqrt{P_{d,v}^2+Q_{d,v}^2}}{V_v^{\mathrm{nom}}},\\
   \alpha_i&=\frac{\operatorname{median}_j(d_j)}{d_i},
 \end{aligned}
\]
where $P_{d,v}$ and $Q_{d,v}$ are active and reactive demands, $V_v^{\mathrm{nom}}$ is nominal voltage in kilovolts, and the median is over subinstances in the same composite.
We multiply demands, shunt injections, generator outputs and bounds, line ratings, and line charging by $\alpha_i$, divide internal resistance and reactance by $\alpha_i$, and reparameterize generator costs.
Voltage data, tap ratios, topology, and coupling branches are unchanged.

\begin{table}[htbp]
  \centering
  \scriptsize
  \caption{Counts of SCOPF runs returning a feasible solution before and after normalization across voltage levels; 100 composites per method; overload slack permitted. Bold: setting used in the main experiment.}
  \label{tab:eval-normalization}
  \setlength{\tabcolsep}{4pt}
  \begin{tabular}{lrrrrrr}
  \toprule
  & R & HD & KK & OF-H & OF-B & Total \\
  \midrule
  \textbf{Without normalization} & 64 & 70 & 53 & 61 & 62 & 310 \\
  With normalization    & 58 & 70 & 46 & 57 & 57 & 288 \\
  \bottomrule
\end{tabular}

\end{table}

Across the 343 extracted subinstances, scale factors range from \num{0.000726} to \num{885.339}.

\FloatBarrier

\end{document}